\RequirePackage{fix-cm}
\documentclass{article}       

\makeatletter
\def\cl@chapter{\@elt {theorem}}
\makeatother
\usepackage[top=1 in, bottom=1 in, left=1 in, right = 1 in]{geometry}
\usepackage{amsmath}
\usepackage{amsbsy}
\usepackage{amssymb}
\usepackage{setspace}
\usepackage{graphicx}
\usepackage{algorithm}
\usepackage{algorithmic}
\usepackage{hyperref}
\usepackage{cite}
\usepackage[numbers]{natbib}
\usepackage[font=scriptsize,labelfont=bf]{caption}
\usepackage[font=scriptsize,labelfont=bf]{subcaption}
\usepackage{textcomp}
\usepackage{algorithm}
\usepackage{algorithmic}
\usepackage{animate}
\usepackage{cleveref,graphicx,subcaption,tikz,wrapfig}
\newcommand{\sign}{\text{sign}}

\newcommand{\abs}[1]{\left| #1  \right|}

\newcommand{\R}{\mathbb{R}}

\newcommand{\bx}{\boldsymbol{x}}

\newcommand{\ba}{\boldsymbol{a}}

\renewcommand{\o}{\omega}

\newcommand{\eps}{\varepsilon}

\newcommand*{\defeq}{\mathrel{\vcenter{\baselineskip0.5ex \lineskiplimit0pt
                     \hbox{\scriptsize.}\hbox{\scriptsize.}}}%
     		     =}

\def\XXint#1#2#3{{\setbox0=\hbox{$#1{#2#3}{\int}$ }
\vcenter{\hbox{$#2#3$ }}\kern-.6\wd0}}

\begin{document}

\title{\Large A Rotating-Grid Upwind Fast Sweeping Scheme for a Class of Hamilton-Jacobi Equations}

\author{Christian Parkinson\thanks{Department of Mathematics, University of Arizona, Tucson, AZ, 85721 (\emph{chparkin@math.arizona.edu}).}}

\date{ }



\maketitle

\begin{abstract}
We present a fast sweeping method for a class of Hamilton-Jacobi equations that arise from time-independent problems in optimal control theory. The basic method in two dimensions uses a four point stencil and is extremely simple to implement. We test our basic method against Eikonal equations in different norms, and then suggest a general method for rotating the grid and using additional approximations to the derivatives in different directions in order to more accurately capture characteristic flow. We display the utility of our method by applying it to relevant problems from engineering. 
\end{abstract}

\section{Introduction} \label{sec:intro}

The general Hamilton-Jacobi (HJ) equation in $d$-dimensions is given by \begin{equation}\label{eq:HJgeneral}
H(x, \nabla\phi(x)) = 0, \,\,\,\, x \in \Omega
\end{equation} where $ \Omega\subset \R^d$ and $H:  \Omega \times \R^d \to \R$ is the \emph{Hamiltonian} function. Along with equation \eqref{eq:HJgeneral}, one is often supplied boundary data $\phi(x) = g(x)$ on a set $\Gamma \subset \R^d$, which typically has dimension smaller than $d$. Common scenarios are $\Gamma = \partial \Omega$ or $\Gamma = \{x_0\}$, a single point. These equations have diverse application in fields including traffic modeling \cite{Luo1}, medical imaging \cite{Luo2}, path-planning \cite{Parkinson}, and dynamic visibility \cite{Kao,ObermanVisibility,Tsai} to name a few. 

The fast sweeping method is a type of finite difference scheme used to approximate \eqref{eq:HJgeneral}.  The basic strategy involves discretizing the domain and devising update rules \begin{equation} 
u_i = F_i(u_j\rvert_{j \in N(i)})
\end{equation} that locally approximate the equation at grid nodes $i$, where $N(i)$ is comprised of the nodes in some neighborhood of node $i$. Using these update rules, one sweeps through the domain in the Gauss-Seidel manner, iteratively updating the solution values at grid nodes until convergence. As far as this author can discern, the fast sweeping method was first used by Bou\'{e} and Dupuis \cite{Boue} and Zhao et al. \cite{ZhaoOsher}. Shortly afterwards, there was much work on developing fast sweeping methods for different types of Hamiltonians, and using different strategies for numerical approximation \cite{Kao20042,Kao2004,Kao2005,TsaiOsherSweep,Zhao}. Subsequent effort was devoted to adapting fast sweeping methods to irregular grids \cite{Qian2,Qian1}, improving the accuracy \cite{Luo3,Luo4}, and extending them to other equations, such as conservation laws \cite{Engquist1,Engquist2}. Luo and Zhao \cite{LuoZhao} provide a nice overview of fast sweeping methods, which we will refer to in \cref{sec:analysis}. 

Besides fast sweeping schemes, other grid-based methods used to approximate steady-state HJ equations can be largely divided into two categories. The first category is fast marching methods for monotonically advancing fronts, pioneered by Tsitsiklis \cite{Tsitsiklis}. These methods---as well as their generalization to ordered upwind methods---rely on a single-pass, Dijkstra-type algorithm to update the solution value at grid nodes as characteristics flow outward from boundary data \cite{IanMitchell, Alton, Sethian1, SethVlad2, SethVlad1}. Besides these, Bornemann and Rasch \cite{BR} proposed a variational method based on the Hopf-Lax formula. Their approach is to localize the HJ equation to finitely many simplices, approximate the solution with linear elements, and solve a discrete version of the Hopf-Lax formula. Their method is similar in spirit to fast marching methods in that it involves updating nodes in a specific order. However, it relies on a Gauss-Seidel iteration, rather than a single pass update. The second category is time-dependent methods. Osher showed that in many cases one can recast the steady-state HJ equation in a time-dependent manner  \cite{Osher1993}. There are very general methods which can approximate time-dependent HJ equations at high accuracy, and also allow for non-monotonic flow of information \cite{JiangWENO,OsherShu,Shu}. More recently, there has been increased interest in algorithms for numerical solutions of HJ equations which break the curse of dimensionality. These typically rely on Hopf-Lax or Lax-Oleinik type formulas for time-dependent HJ equations, and use optimization routines to approximate the solution at individual points \cite{YTChow,Darbon,LinCurse}. However, due to the wide applicability and relative ease of both implementation and analysis, fast sweeping methods have remained a popular option for approximating solutions of steady-state HJ equations. 

We present an exceedingly simple fast sweeping scheme for a class of Hamilton-Jacobi equations arising from optimal control theory. For simplicity of exposition, we develop our method in two spatial dimensions. The method applies in higher dimensions, though for dimensions $d > 3$, one will encounter the curse of dimensionality. In two dimensions, our most basic method includes a four-point stencil on a rectangular grid, using only the ordinary forward and backward difference operators. We then describe a general method for using rotated coordinates to improve the accuracy of the scheme. We implement our method with special application toward Eikonal equations in different norms, and also mention a few other applications. Because one of the strengths of our method is ease of implementation, we compare it with the Lax-Friedrichs sweeping scheme \cite{Kao2004}, another easily implementable method. \\

\section{Hamilton-Jacobi Equations in Optimal Control Theory}\label{sec:HJOpt} 

We will address a specific class of Hamilton-Jacobi equations arising from deterministic optimal control theory. A basic problem in optimal control theory is to choose the best control plan $\ba: [0,T] \to A$ to steer a trajectory $\bx$ obeying \begin{equation} \label{eq:motion} \begin{split}
&\dot \bx(t) = f(\bx(t),\ba(t),t), \,\,\,\, 0 < t \le T, \\
&\bx(0) = x_0,
\end{split}
\end{equation} to an optimal destination $\bx(T)$. Here $A\subset \R^m$ is the set of admissible control actions and $f: \R^d \times \R^m \times [0,T]\to \R^d$ is a function describing the dynamics along the trajectory. The ``optimal destination" is determined in view of a cost functional \begin{equation} 
C[\bx(\cdot),\ba(\cdot)] = g(\bx(T)) + \int^T_0 r(\bx(t),\ba(t),t)dt 
\end{equation} that one wishes to minimize. The function $r : \R^d \times \R^m \times [0,T]\to \R$ accounts for a running cost along the trajectory, and $g: \R^d\to\R$ is the exit cost. While it is not necessary in all cases, we will assume that $r,g \ge 0$, which is common in many applications where cost cannot be negative. To analyze this problem using dynamic programming \cite{Bellman,Bellman2}, one defines the value function $\phi: \R^d \times [0,T] \to\R$ by \begin{equation} \label{eq:value}
\phi(x,t) \defeq \inf_{\bx(\cdot),\ba(\cdot)} C_{x,t} [\bx(\cdot),\ba(\cdot)]
\end{equation} where $C_{x,t}[\bx(\cdot),\ba(\cdot)]$ is the remaining cost functional, restricted to trajectories $\bx$ on the time interval $(t,T]$ and satisfying $\bx(t) = x$. Thus $\phi$ is the optimal remaining cost for a trajectory that is at position $x$ at time $t$. Under mild conditions on the data, this value function is the unique viscosity solution \cite{CrandallLions} of the terminal value Hamilton-Jacobi-Bellman equation \cite{BardiCapuzzo,Barles} \begin{equation} \label{eq:HJBtimeDep}
\begin{split}
&\phi_t(x,t) + \inf_{a \in A} \Big\{\langle f(x,a,t), \nabla \phi(x,t)\rangle + r(x,a,t) \Big\} = 0, \\
&\phi(x,T) = g(x).
\end{split}
\end{equation} Note that the viscosity solution of \eqref{eq:HJBtimeDep} should remain non-negative: by \eqref{eq:value}, $\phi$ is non-negative whenever $r$ and $g$ are non-negative. 

We observe that \eqref{eq:HJBtimeDep} is of the form \eqref{eq:HJgeneral} if we consider generalized coordinates $\tilde x = (t,x)$ and $\nabla_{\tilde x} = (\partial_t, \nabla_x)$. In this case $\Omega = \R^d \times [0,T)$ and $\Gamma = \R^d \times \{T\}$. Thus this can be analyzed in the framework of the more general equation \eqref{eq:HJgeneral}, but time-dependent equations like \eqref{eq:HJBtimeDep} are so ubiquitous in application that they are often analyzed independently. Indeed, in their two original papers, Crandall and Lions established the notion of viscosity solutions specifically for time-dependent Hamilton-Jacobi equations \cite{CrandallLions,CrandallLions2}, and later the theory was extended to more general equations; see, for example, \cite{CrandallIshiiLions1992}. 

\subsection{Our Class of Equations} \label{sec:ours} 

We restrict our focus to a special class of optimal control problems. We consider the case that the dynamic function $f$ does not depend explicitly on $t$, and the running cost function $r$ does not depend explicitly on either $t$ or $\ba(\cdot)$. The removal of the explicit dependence on $t$ is not a particularly stringent condition; this is very natural in many applications. Removing the dependence of $r$ on $\ba(\cdot)$ is a more serious restriction. For example, this will exclude essentially any problem from mathematical finance where the control variable could represent the fraction of capital one wishes to invest or the amount of goods a company would like to produce \cite{Pham}. In this case, the cost and profit very explicitly depend on the value of the control variable. However, control problems of our type still have diverse application. Minimal-time path-planning \cite{Parkinson} and reach avoid games \cite{TomlinReachAvoid} are two classical problems in applied optimal control theory that fit into this framework. Otherwise, four of the five examples given by Evans \cite[chap.~1]{EvansControlNotes} fall into this category. This includes the moon lander problem, optimally stopping a pendulum, and a model for growth of ant colonies originally proposed by Oster and Wilson \cite{OsterWilson}.

When neither $f$ nor $r$ depend on $t$, one can neglect the time horizon $T$ and formulate a steady-state Hamilton-Jacobi-Bellman equation for the value function. Given that $r$ does not depend on $\ba(\cdot)$, this takes the form \begin{equation}  \label{eq:ourEquation1}
-r(x) = \inf_{a\in A} \left\{ \langle f(x,a),\nabla \phi(x)\rangle\right\},
\end{equation} or alternately \begin{equation}  \label{eq:ourEquation2}
-r(x) = \inf_{a\in A} \left\{ \sum^d_{\ell=1} f_\ell(x,a) \phi_{x_\ell}(x) \right\} 
\end{equation} where $x = (x_1,\ldots, x_d)$ and $f(x,a) = (f_1(x,a),\ldots,f_d(x,a))$. We focus on numerical solutions for this equation with boundary data $\phi(x) = g(x)$ on a set $\Gamma \subset \R^d$. For example, in the case of optimal-time path-planning, we will take $\Gamma = \{x_f\}$, where $x_f \in \R^d$ is the desired ending point, and let $\phi(x_f) = 0$. This signifies that paths ending at the desired location incur no exit cost, while other paths are not admissible (i.e., they incur infinite cost).

Many classical Hamilton-Jacobi equations can be expressed in this form. Notably, the Eikonal equation \begin{equation} \label{eq:eikonal}
1 = v(x) \abs{\nabla \phi(x)}
\end{equation} is of this form. The travel-time function for isotropic motion $\dot \bx(t) = v(\bx(t))\ba(t)$, where $\ba(\cdot)$ is a unit vector, is the viscosity solution of this equation, and in the case that $v(x) \equiv 1$, this yields a signed distance function \cite{OsherFedkiw}. Assuming $v > 0$, equation \eqref{eq:eikonal} can be re-written \begin{equation}
-1/v(x) = \inf_{a \in \mathbb S^{d-1}} \Big\{ a \cdot \nabla \phi  \Big\}
\end{equation} whereupon casting the equation in the form \eqref{eq:ourEquation2} is accomplished by parameterizing the unit sphere $\mathbb S^{d-1}$. For example in dimension $d =2$, we have \begin{equation}\label{eq:eikonal2d}
-1/v(x,y) = \inf_{a \in [0,2\pi)} \Big\{ \phi_x \cos(a) + \phi_y\sin(a)\Big\},
\end{equation} or in dimension $d = 3$, \begin{equation} \label{eq:eikonal3d}
-1/v(x,y,z) = \inf_{a,b} \Big\{\phi_x\cos(a)\cos(b) + \phi_y \sin(a)\cos(b) + \phi_z\sin(b)  \Big\},
\end{equation} where $(a,b)\in[0,2\pi)\times [-\pi/2,\pi/2]$ represent the $xy$-planar angle and the angle of inclination from the $xy$-plane, respectively. We return to Eikonal equations when testing our method in \cref{sec:eikonalBasic} and \cref{sec:eikonalRot}.

\section{A Basic Fast Sweeping Scheme for \eqref{eq:ourEquation2}} \label{sec:ourBasicScheme} As stated in \cref{sec:intro}, for simplicity of exposition, we will describe our fast sweeping scheme in dimension $d = 2$. We consider a rectangular domain $[x_{\text{min}}, x_{\text{max}}] \times [y_{\text{min}}, y_{\text{max}}]$ and a uniform grid discritization with $I+1$ points in the $x$-direction, and $J+1$ points in the $y$-direction. Thus the grid is given by \begin{equation}\label{eq:grid}\begin{split}
x_i &\defeq x_{\text{min}} + i\Delta x, \,\,\,\,\, \Delta x = \frac{x_{\text{max}} - x_{\text{min}}}{I}, \,\,\,\, i = 0,1,\ldots,I,\\
y_j &\defeq y_{\text{min}} + j\Delta y, \,\,\,\,\, \Delta y = \frac{y_{\text{max}} - y_{\text{min}}}{J}, \,\,\,\, j = 0,1,\ldots,J.
\end{split}
\end{equation} In two-dimensions, the equation of interest is \begin{equation} \label{eq:ourEquation2d}
-r(x,y) = \inf_{a \in A} \Big\{ f_1(x,y,a) \phi_x(x,y) + f_2(x,y,a)\phi_y(x,y) \Big\}.
\end{equation}  Let $\phi_{ij}$ be the numerical approximation to $\phi(x_i,y_j)$, and for a fixed $a \in A$, let $f_{\ell,ij}(a) = f_\ell(x_i,y_j,a)$ for $\ell = 1,2$. Further let \begin{equation}\label{eq:signs}
\xi_{\ell,ij}(a) = \text{sign}(f_\ell(x_i,y_j,a)), \,\,\,\,\,\,\, \ell = 1,2. 
\end{equation} Then the upwind approximations to the derivatives are given by \begin{equation}\label{eq:upwindApprox}\begin{split}
\Big(f_{1}(x,y,a)\phi_x(x,y)\Big)_{ij} &= \abs{f_{1,ij}(a)} \frac{\phi_{i+\xi_{1,ij}(a),j} - \phi_{ij}}{\Delta x},\\ 
\Big(f_{2}(x,y,a)\phi_y(x,y)\Big)_{ij} &= \abs{f_{2,ij}(a)} \frac{\phi_{i,j+\xi_{2,ij}(a)} - \phi_{ij}}{\Delta y}.
\end{split}
\end{equation} Supposing that $a$ is the correct control value at the node $(i,j)$, we can insert these approximations into \eqref{eq:ourEquation2d} to arrive at \begin{equation}\label{eq:discEq1}
-r_{ij} = \abs{f_{1,ij}(a)} \frac{\phi_{i+\xi_{1,ij}(a),j} - \phi_{ij}}{\Delta x} + \abs{f_{2,ij}(a)} \frac{\phi_{i,j+\xi_{2,ij}(a)} - \phi_{ij}}{\Delta y},
\end{equation} where $r_{ij} = r(x_i,y_j)$. Isolating $\phi_{ij}$, we see that \begin{equation}\label{eq:updateRule}
\phi^*_{ij}(a) = \frac{r_{ij} + \frac{\abs{f_{1,ij}(a)}}{\Delta x}\phi_{i+\xi_{1,ij}(a),j} +  \frac{\abs{f_{2,ij}(a)}}{\Delta y}\phi_{i,j+\xi_{2,ij}(a)}}{\frac{\abs{f_{1,ij}(a)}}{\Delta x} + \frac {\abs{f_{2,ij}(a)}}{\Delta y}}
\end{equation} is a first-order upwind approximation to equation \eqref{eq:ourEquation2d}, when $a$ is the correct control value at node $(i,j)$. This suggests the fast sweeping scheme detailed in \cref{alg:sweepingScheme}.

\begin{algorithm}[t!]
    \caption{A fast sweeping scheme to solve \eqref{eq:ourEquation2d}}
    \label{alg:sweepingScheme}
    \begin{algorithmic}
    \STATE {\bf Initialization:} Input boundary data (a function $g$ and set $\Gamma$), a grid discretization as in \eqref{eq:grid}, and a small error tolerance $\eps > 0$. Initialize $\phi^0_{ij} = g(x_i,y_j)$ for the grid nodes corresponding to $\Gamma$ and $\phi^0_{ij} = +\infty$ (or some large positive number) for all other grid nodes. Initialize $\phi^1_{ij} = 0$ at all grid points, and $n = 1$. \\
    \STATE
        \WHILE{$\|\phi^n - \phi^{n-1}\| > \eps$}
        \STATE
        \STATE Assign $\phi^n_{ij} \leftarrow \phi^{n-1}_{ij}$ for all $(i,j)$. 
        \STATE
        \FOR{$i = 1$ \TO  $I-1$}
        \FOR{$j = 1$ \TO  $J-1$}
        \STATE For each $a \in A$, compute $$\phi^*_{ij}(a) \leftarrow \frac{r_{ij} + \frac{\abs{f_{1,ij}(a)}}{\Delta x}\phi^n_{i+\xi_{1,ij}(a),j} +  \frac{\abs{f_{2,ij}(a)}}{\Delta y}\phi^n_{i,j+\xi_{2,ij}(a)}}{\frac{\abs{f_{1,ij}(a)}}{\Delta x} + \frac {\abs{f_{2,ij}(a)}}{\Delta y}}.$$ 
        \STATE
        \STATE Assign $\phi^n_{ij} \leftarrow \min\{\min_{a} \phi^*_{ij}(a), \phi^{n-1}_{ij} \}$
        \ENDFOR
        \ENDFOR
        \STATE
        \STATE Repeat the above {\bf for} loops, sweeping in alternating directions until all combinations of sweeping directions have been completed (a total of 4 sweeps). 
        \STATE 
        \STATE Assign $n \leftarrow n+1$
        \STATE
        \ENDWHILE
	\STATE
	\STATE {\bf return} the values $\phi^{\text{end}}_{ij}$ for all $(i,j)$
    \end{algorithmic}
\end{algorithm}

We include some comments regarding the algorithm. First, at each iteration, we sweep through the indices in alternating directions until all combinations of sweeping directions have been performed. Thus each iteration consists of four sweeps; in MATLAB notation: \begin{enumerate}
\item[(1)] $i = 1\,:\, I-1, \hphantom{ -1 \,:\,} \,\,\,\,\,\, j = 1 \,:\, J-1,$
\item[(2)] $i = 1\,:\, I-1, \hphantom{ -1 \,:\,} \,\,\,\,\,\, j = J-1 \,:\, -1 \,:\, 1,$
\item[(3)] $i = I-1\,:\,-1 \,:\, 1, \,\,\,\,\, j = J-1 \,:\, -1\,:\,1,$
\item[(4)] $i = I-1\,:\, -1\,:\,1, \,\,\,\,\, j = 1 \,:\, J-1.$
\end{enumerate} Generally, in dimension $d$, there will be $2^d$ sweeps in each iteration. Second, it is important  that we assign $\phi^n_{ij} \leftarrow \phi^{n-1}_{ij}$ at the beginning of each iteration and then operate only with $\phi^n_{ij}$. This ensures that sweeping is carried out in the Gauss-Seidel sense: updating values, and then using the most recently updated values to resolve the ensuing values. Third, for the convergence criterion, we use the $L^\infty$-norm so that the iteration halts when $\|\phi^n-\phi^{n-1}\| = \max_{ij}\abs{\phi_{ij}^n-\phi_{ij}^{n-1}} \le \eps$ for some prescribed tolerance $\eps$, though other criteria could be used. Fourth, the scheme is fully upwind meaning that numerical characteristics flow away from the boundary set $\Gamma$. If $\Gamma$ corresponds to the computational boundary, then information flows into the domain. If $\Gamma$ is contained in the computational domain, then characteristics will flow out of the computational boundary. In this case, no special considerations are necessary at the computational boundaries. The values at the boundary nodes will remain large, but will not affect the solution at interior nodes. In this way, our scheme is similar to Godunov-inspired methods such as \cite{TsaiOsherSweep}. In a different approach, Kao et al. \cite{Kao2004} devise a sweeping method with a Lax-Friedrichs Hamiltonian, wherein added numerical diffusion will cause information to seep into the domain from the computational boundary, requiring special consideration. We will discuss the Lax-Friedrichs sweeping method in more detail later. 

Perhaps the most important notes regard the minimization over $a \in A$, which takes place at each grid point in each sweep. A single iteration requires this minimization to be resolved roughly $4IJ$ times. Because of this, the shape of $A$ is somewhat crucial to the algorithm. For example, in the Eikonal equation, we have $A = \mathbb S^1$, meaning this optimization is performed over a continuous set. One can either discretize the set and choose from finitely many values, or introduce an optimization routine of their choosing. Either way, this is likely to represent the largest computational burden. The algorithm performs extraordinarily well when $A$ is finite. For example, this occurs in bang-bang control problems, where the optimal controls switch between finitely many control values \cite{BangBang}. One application of this is in kinematic models for simple self-driving cars \cite{Dubins,ReedsShepp}. Takei and Tsai were the first to analyze this problem in the Hamilton-Jacobi setting \cite{TakeiTsai2,TakeiTsai1}, and they used a sweeping scheme just like ours. We will return to the example of self-driving cars in \cref{sec:otherApps}.
  
\subsection{Upwinding, Monotonicity \& Convergence}\label{sec:analysis} Luo and Zhao \cite{LuoZhao} discuss and analyze fast sweeping methods in some generality. In particular, they consider \eqref{eq:HJgeneral} with a Hamiltonian $H$ that is \begin{itemize} \item[$(i)$] continuous on $\Omega \times \R^n$, \item[$(ii)$] convex and coercive in $\nabla \phi$, \item[$(iii)$] compatible, in that $H(x,0) \le 0$ for $x \in \Omega$. \end{itemize} Under these conditions and some mild conditions on the boundary data $g$, they prove that if a fast sweeping scheme is consistent, monotone, and obeys a causality condition, then the approximate solution produced by the scheme will converge to the viscosity solution of the Hamilton-Jacobi equation under grid refinement.  

\begin{figure}[t!]
\centering
\includegraphics[width=0.6\textwidth,trim=50 50 50 50,clip]{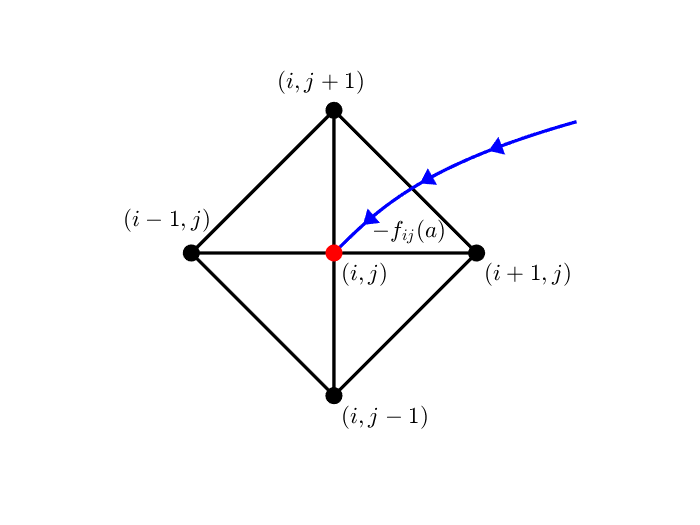}
\caption{The causality condition specifies that the nodes used to approximate $(\nabla \phi)_{ij}$ form a polygon containing the characteristic (blue) flowing into $(i,j)$. Here, one would use nodes $(i,j), \,\, (i+1,j), \,\, (i,j+1)$. The characteristic direction is given by $-f_{ij}(a)$ if $a$ is the correct control value at the grid node.}
\label{fig:upwind}
\end{figure}

An annoying but necessary facet of the theory of viscosity solutions is that orientation matters. Formally, the viscosity solution of $H(x,\nabla \phi(x))=0$ is the negative of the viscosity of $-H(x,\nabla\phi(x))=0$. Our orientation is reversed from that in \cite{LuoZhao} but modulo some sign changes and inequality flips, the analysis is the same. Our scheme is consistent to first order, as can be shown by a simple Taylor expansion. In our case, the monotonicity requirement is trivially satisfied since the update rule \eqref{eq:updateRule} is clearly non-decreasing in the values at the surrounding grid nodes. The causality condition states in essence that the characteristic flowing into grid node $(i,j)$ is contained in the polygon formed by the nodes used for the finite difference approximations at $(i,j)$. This is illustrated in \cref{fig:upwind}, where the characteristic curve (blue) enters from the positive-$x$ and positive-$y$ direction, specifying that one should use nodes $(i,j), \,\, (i+1,j), \,\, (i,j+1)$ to approximate $(\nabla \phi)_{ij}$. For us, the causality condition corresponds exactly to the upwind approximations \eqref{eq:upwindApprox}. Note that because of the negative sign in the equation, the characteristic direction at $(x,y)$ is  $-f(x,y,a)$ when $a$ is the correct control value at $(x,y)$. Thus our scheme fits into their framework, and we have convergence to the viscosity solution of \eqref{eq:ourEquation1} as the grid parameters go to zero.

Determining the order of covergence is subtle. Classical proofs of convergence for numerical solutions of Hamilton-Jacobi equations depend not only on the order of local truncation error, but also on the regularity of the viscosity solution \cite{BarlesJakobsen,BarlesSouganidis,Souganidis1985}. Typically one can guarantee convergence at order no less than $1/2$ when the scheme is consistent at order $1$. However, one often sees full first-order convergence in regions where the solution is smooth \cite{LuoZhao}, and in some cases, one can achieve higher order accuracy using techniques such as ENO or WENO schemes \cite{JiangWENO,OsherShu,Shu, ZhangWENOFSM}, though the application of these concepts to fast sweeping methods presents some challenges. We discuss this further in \cref{sec:WENO}.  \\

\subsection{Application of the Basic Method to Eikonal Equations} \label{sec:eikonalBasic} To empirically study error and convergence, we test our method on three different Eikonal equations: \begin{equation}\label{eq:eikonalpnorm}
1 = \|\nabla \phi(x) \|_{p}
\end{equation} where $p = 1,2, \infty$. Given the boundary data $\phi(0) = 0$, we see that the unique (positive) viscosity solution of \eqref{eq:eikonalpnorm} is $\phi_p(x) = \|x\|_{p'}$ where $\frac 1 p + \frac 1{p'} = 1$. This fact can be intuited from the ensuing optimal control problem, and essentially follows from the dual definition of the norm: \begin{equation} \label{eq:dualNorm}
\|z\|_p = \sup_{\|a\|_{p'}\le 1 } \langle z,a\rangle.
\end{equation} However, proving this in full generality is surprisingly intricate. A discussion of such equations is included in \cite{ObermanVladimirsky}, and a full analysis is given in \cite{CaffarelliCrandall}. 

Each of these equations is solved by travel time function for a minimal-time path-planning problem of the form above. Indeed, consider the equation of motion \begin{equation} \label{eq:motionEikonal} 
\dot \bx(t) = \ba(t), \,\,\,\,\,\, \ba(\cdot) \in B^{(p')}_1,
\end{equation} where $B_1^{(p')}$ is the unit ball in the $p'$-norm (centered at the origin). If we pair this equation with the cost functional \begin{equation}
C[\bx(\cdot),\ba(\cdot)] = \iota_0(\bx(T)) + \int^T_0 1 \, dt
\end{equation} where $\iota_0$ is the convex indicator of the origin ($0$ at the origin; $+\infty$ elsewhere) and allow for infinite horizon time, then the Hamilton-Jacobi-Bellman equation for the value function is the $p$-norm Eikonal equation \eqref{eq:eikonalpnorm}, and the optimal control plan steers the trajectory to the origin in the minimal possible time, where distance from the origin is computed in the $p'$-norm. In particular, since the unit ball has finitely many extreme points in the case that $p'=\infty$ or $p'= 1$, this leads to a bang-bang control problem for $p = 1$ or $p = \infty$. 

In two-dimensions, equation \eqref{eq:eikonal2d} shows that the $2$-norm Eikonal equation can be written in the form \eqref{eq:ourEquation2d}.  We can write the other equations in this form as well. For $p =1$, we have \begin{equation}\label{eq:eikonal1norm}
 -1 = \inf_{a_1,a_2 \in \{\pm 1\}} \Big\{ a_1 \phi_x(x,y) + a_2 \phi_y(x,y) \Big\}
\end{equation} and for $p = \infty$, we have \begin{equation}\label{eq:eikonalInfnorm}
 -1 = \inf_{a \in \{\pm e_1,\pm e_2\}} \Big\{ a_1 \phi_x(x,y) + a_2 \phi_y(x,y) \Big\},
\end{equation} where in the latter equation, $e_1,e_2$ are the standard basis vectors, and $a = (a_1,a_2)$.

We would like derive the specific update formula \eqref{eq:updateRule} for each of these cases. For the ordinary Eikonal equation in the $2$-norm, we find \begin{equation}\label{eq:updateRule2norm}
\phi^{*,2}_{ij}(a) = \frac{1 + \frac{\abs{\cos(a)}}{\Delta x}\phi^{n,2}_{i+\text{sign}(\cos(a)),j} + \frac{\abs{\sin(a)}}{\Delta y} \phi^{n,2}_{i,j+\text{sign}(\sin(a))}}{\frac{\abs{\cos(a)}}{\Delta x} + \frac{\abs{\sin(a)}}{\Delta y}}.
\end{equation} and use the update $\phi^{n,2}_{ij} = \min\{ \min_{a \in [0,2\pi)} \phi^{*,2}_{ij}(a), \phi^{n-1,2}_{ij}\}$. To use this update, we will need to resolve the minimization over $a \in [0,2\pi)$. To do so, we simply sample $a = 2\pi k /K$ for $k=0,\ldots, K-1$ and choose the minimum from these finitely many points. In our tests, we fix $K = 400$. This will incur some small error. We discuss this briefly below.

For the $1$-norm and $\infty$-norm equations, we can explictly write the update rule by considering all possible combinations of control variables. For the case $p=1$, we have \begin{equation} \label{eq:updateRule1norm} \begin{split}
\phi^{n,1}_{ij} = \min\bigg\{  \phi^{n-1,1}_{ij}, &\frac{1 + \frac{1}{\Delta x} \phi^{n,1}_{i+1,j} + \frac{1}{\Delta y} \phi^{n,1}_{i,j+1}}{\frac{1}{\Delta x} + \frac 1 {\Delta y}},  \frac{1 + \frac{1}{\Delta x} \phi^{n,1}_{i-1,j} + \frac{1}{\Delta y} \phi^{n,1}_{i,j+1}}{\frac{1}{\Delta x} + \frac 1 {\Delta y}},\\
& \frac{1 + \frac{1}{\Delta x} \phi^{n,1}_{i+1,j} + \frac{1}{\Delta y} \phi^{n,1}_{i,j-1}}{\frac{1}{\Delta x} + \frac 1 {\Delta y}}, \frac{1 + \frac{1}{\Delta x} \phi^{n,1}_{i-1,j} + \frac{1}{\Delta y} \phi^{n,1}_{i,j-1}}{\frac{1}{\Delta x} + \frac 1 {\Delta y}} \bigg\}.
\end{split} 
\end{equation} 

In the $p = \infty$ case, the update is even simpler since one of $a_1,a_2$ in \eqref{eq:eikonalInfnorm} is zero. Plugging the values into the general update formula \eqref{eq:updateRule} and clearing the denominator yields \begin{equation}\label{eq:updateRuleInfnorm}
\phi^{n,\infty}_{ij} = \min\Big\{\phi^{n-1,\infty}_{ij},\,\, \Delta x + \phi^{n,\infty}_{i+1,j},\,\, \Delta x + \phi^{n,\infty}_{i-1,j},\,\, \Delta y + \phi^{n,\infty}_{i,j+1},\,\, \Delta y + \phi^{n,\infty}_{i,j-1}  \Big\}.
\end{equation} We note that \eqref{eq:updateRuleInfnorm} is perfectly satisfied by the exact solution $\phi_\infty(x,y) = \|(x,y)\|_1 = \abs{x} + \abs y$, and thus when $p=\infty$, our scheme will solve the equation exactly, so long as the origin is a grid node. Otherwise, the error in the approximation will only depend on the distance from the origin to the nearest grid node in each direction.

\begin{figure}[t!]
\begin{subfigure}{0.3\textwidth}
\centering 
\includegraphics[width=\textwidth,trim=43 15 35 5,clip]{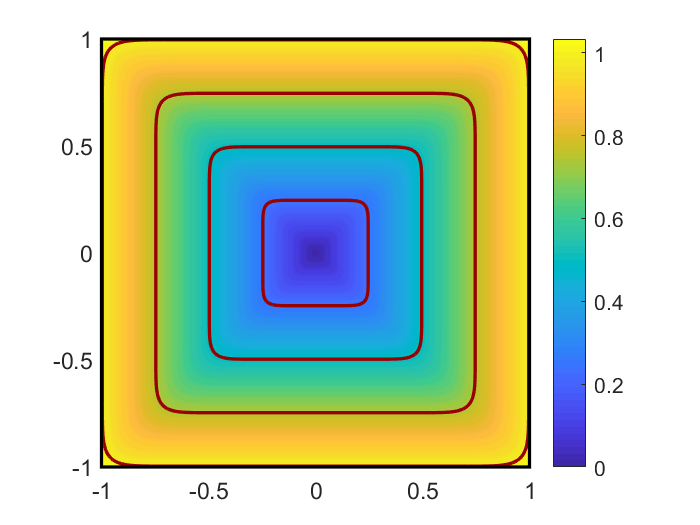}
\caption{Approx. soln., $p=1.$}
\label{fig:OneNormSoln}
\end{subfigure}~
\begin{subfigure}{0.3\textwidth}
\centering 
\includegraphics[width=\textwidth,trim=43 15 35 5,clip]{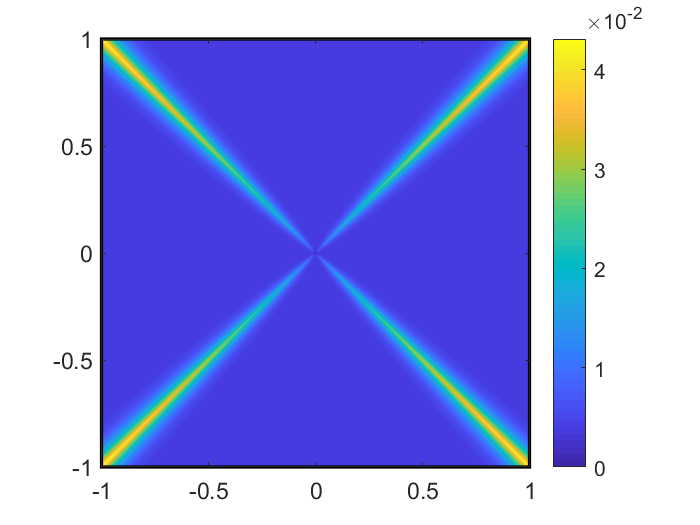}
\caption{Error when $p=1$.}
\label{fig:OneNormErr}
\end{subfigure}~
\begin{subfigure}{0.36\textwidth}
\centering
\begin{tabular}{r r r}
$I,J$ & $L^\infty$ Err. & Conv.  \\ [0.5ex]
\hline \hline 
  50& 	 1.4057e-01& 	 ---\\
 100&	 9.3988e-02& 	 0.5807\\
 200& 	 6.3636e-02& 	 0.5626\\
 400& 	 4.3544e-02& 	 0.5474\\
 800& 	 3.0049e-02& 	 0.5352\\
1600& 	 2.0872e-02& 	 0.5257\\ [0.5ex]
\hline 
\end{tabular}
\caption{Conv. table for $p=1$.}
\label{fig:OneNormConv}
\end{subfigure} \\
\begin{subfigure}{0.3\textwidth}
\centering 
\includegraphics[width=\textwidth,trim=43 15 35 5,clip]{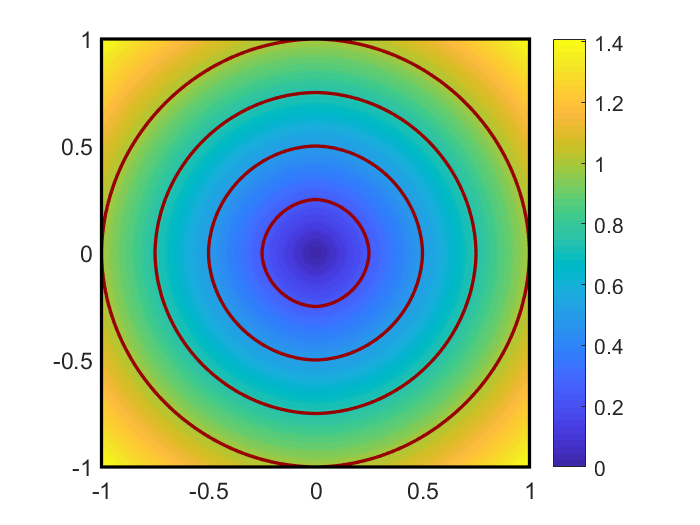}
\caption{Approx. soln., $p=2.$}
\label{fig:TwoNormSoln}
\end{subfigure}~
\begin{subfigure}{0.3\textwidth}
\centering 
\includegraphics[width=\textwidth,trim=43 15 35 5,clip]{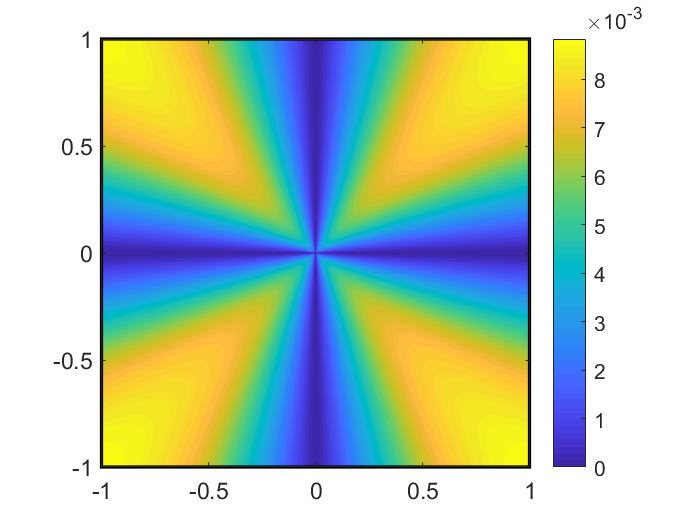}
\caption{Error when $p=2$.}
\label{fig:TwoNormErr}
\end{subfigure}~
\begin{subfigure}{0.36\textwidth}
\centering
\begin{tabular}{r r r}
$I,J$ & $L^\infty$ Err. & Conv.  \\ [0.5ex]
\hline \hline 
  50& 	 4.3754e-02& 	 ---\\
 100& 	 2.6310e-02& 	 0.7338\\
 200& 	 1.5464e-02& 	 0.7666\\
 400& 	 8.9201e-03& 	 0.7938\\
 800& 	 5.0668e-03& 	 0.8160\\
1600& 	 2.8431e-03& 	 0.8336\\ [0.5ex]
\hline 
\end{tabular}
\caption{Conv. table for $p=2$.}
\label{fig:TwoNormConv}
\end{subfigure} \\
\begin{subfigure}{0.3\textwidth}
\centering 
\includegraphics[width=\textwidth,trim=43 15 35 5,clip]{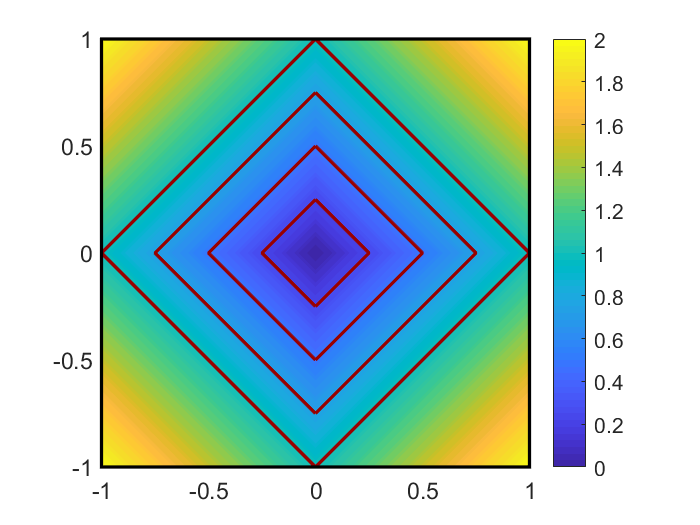}
\caption{Approx. soln., $p=\infty.$}
\label{fig:InfNormSoln}
\end{subfigure}~
\begin{subfigure}{0.3\textwidth}
\centering 
\includegraphics[width=\textwidth,trim=43 15 35 5,clip]{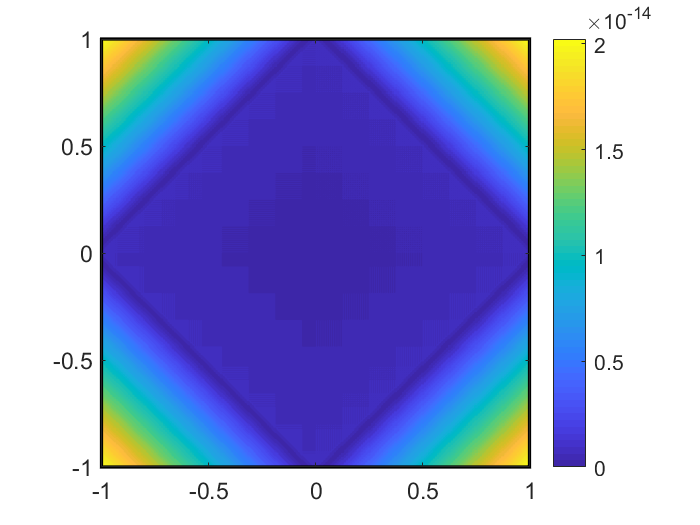}
\caption{Error when $p=\infty$.}
\label{fig:InfNormErr}
\end{subfigure}~
\begin{subfigure}{0.36\textwidth}
\centering
\begin{tabular}{r r r}
$I,J$ & $L^\infty$ Err. & Conv.  \\ [0.5ex]
\hline \hline 
  50& 	 1.7764e-15& 	 ---\\
 100& 	 1.7764e-15& 	 0.0000\\
 200& 	 1.7764e-15& 	 0.0000\\
 400& 	 2.0428e-14& 	 -3.5236\\
 800& 	 4.2633e-14& 	 -1.0614\\
1600& 	 4.2633e-14& 	 0.0000\\ [0.5ex]
\hline 
\end{tabular}
\caption{Conv. table for $p=\infty$.}
\label{fig:InfNormConv}
\end{subfigure} \\
\caption{Approximation of $\| \nabla \phi\|_p = 1$ using our fast sweeping method. Plots display results from the $401\times 401$ grid. Red lines are level sets of the solution.}
\label{fig:eikonalDiffNorms}
\end{figure}

Using these update rules, and the boundary condition $\phi(0,0)=0$, we simulated equation \eqref{eq:eikonalpnorm} for $p=1,2,\infty$. The results are included in \cref{fig:eikonalDiffNorms}. Specifically, results for $p=1$ are included in figures \ref{fig:OneNormSoln}, \ref{fig:OneNormErr}, \ref{fig:OneNormConv}; $p=2$ in figures \ref{fig:TwoNormSoln}, \ref{fig:TwoNormErr}, \ref{fig:TwoNormConv}; and $p=\infty$ in figures \ref{fig:InfNormSoln}, \ref{fig:InfNormErr}, \ref{fig:InfNormConv}. Recall again the exact solution $\phi_p(x,y) = \|(x,y)\|_{p'}.$ The left most figure in each column shows contour plots of the approximate solutions $[-1,1]\times[-1,1]$ with a $401\times 401$ grid, along with level sets of the approximate solutions. The middle figure in each column shows a contour plot of the error in the approximation. The right most figure includes the convergence table in each case. We note that there is a different scale in each plot. 

When $p=1$, the level sets should be perfect squares since these are balls in the $\infty$-norm. At the corners of those squares, the ordinary forward and backward difference operators cannot capture the sharp edges, which leads to some rounding off. Because of this, the error is large along the lines $y=\pm x$, and the order of convergence is roughly $1/2$; the minimal convergence rate guaranteed by the classical theory \cite{BarlesSouganidis,Souganidis1985}.

When $p=2$, the maximum error is less than in the $p=1$ case, and the error itself is more evenly spread throughout the entirety of each quadrant, rather than being focused along specific lines. The convergence rate here is roughly $3/4$, showing improved convergence behavior compared with the $p=1$ case.  An interesting note here is that along the lines $x = 0$ and $y = 0$, the error is effectively zero. This is because the finite difference approximations are focused in those directions, and the cross sections of the exact solution in those directions are linear rays increasing outward from the origin. Thus, for example, when $x>0$, the exact solution satisfies $\phi_2(x+\Delta x,0) = \Delta x + \phi_2(x,0)$, and our discretization captures this relationship with no error. We will return to this line of thought momentarily. Before doing so, we make a further remark regarding the discretization of the control set. Recall, the update rule for the 2-norm Eikonal equation requires that we resolve a minimization problem over $[0,2\pi)$, and to do so we simply discretized the interval into $K=400$ points and chose the minimum from the discrete set. We found empirically that error produced by approximating the control set is smaller than the error in the discrete derivative approximations. To test this, we instead resolved the minimization to a tolerance of $10^{-10}$ using built-in optimization routines in MATLAB. For a $400 \times 400$ grid, the approximate solution found using the exact minimization differed from that found using discrete minimization with $K=400$ by only $1.4 \times 10^{-5}$, whereas the error between the exact solution and each of the approximate solutions was roughly $8.9 \times 10^{-3}$. It bears mentioning that when finding the exact minimum at every point, the algorithm required roughly 150 times the CPU time to resolve the solution. In general, as long as the minimization problem is solved so that the approximation using the exact minimum and the approximation using an approximate minimum differ by no more than $O(\Delta x, \Delta y)$, then the approximation of the control set will not ruin convergence. Beyond that, one must choose how to balance accuracy and efficiency, as well as ease of implementation. To this last point, one of the strengths of this method \emph{is} the ease of implementation, which is why it is particularly suited to problems where the minimization can be resolved explicitly (for example, bang-bang problems such as the 1-norm or $\infty$-norm Eikonal equation or the kinematics of the self-driving car presented in \cref{sec:otherApps}). 

When $p=\infty$, we noted earlier that our scheme should be exact. Indeed, we see that the level sets of the approximate solution are sharp-edged diamonds, exactly mirroring the level sets of $\phi_{\infty}(x,y) = \abs{x} + \abs y$. In this case, the error is near machine-$\varepsilon$, and thus the convergence table is not informative.

\begin{figure}[t!]
\centering 
\includegraphics[width=0.6\textwidth,trim= 0 0 0 0,clip]{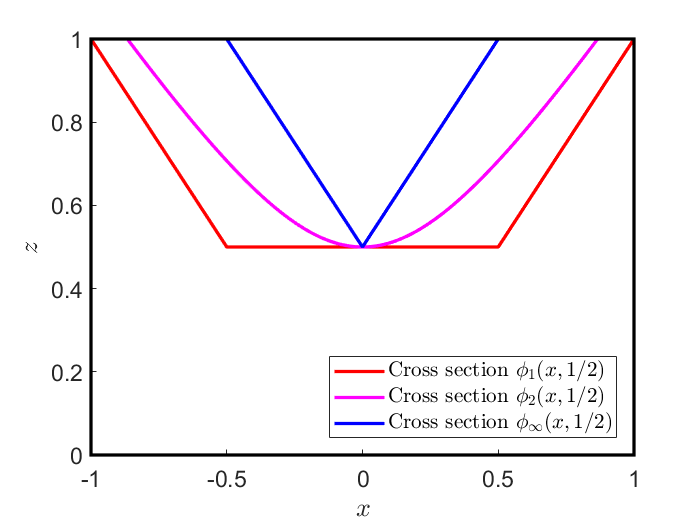}
\caption{Horizontal cross section of $\phi_1,\phi_2,\phi_\infty$ at $y=1/2$.}
\label{fig:crossSec}
\end{figure}

We remarked about the low error along the lines $x=0$ and $y=0$ in the $p=2$ case, and the relationship between this low error and the cross sections of the exact solution along those lines. This remark very closely relates to the improved order of convergence for larger $p$. As $p$ increases (and thus $p'$ decreases), the cross sections of the exact solution $\phi_p(x) = \|x\|_{p'}$ in the vertical or horizontal directions more closely resemble the absolute value function, and thus can be captured more accurately by the finite difference approximations. This is seen in \cref{fig:crossSec}, where we have plotted horizontal cross sections of $\phi_1,\phi_2$ and $\phi_\infty$ at level $y=1/2$. For $\phi_\infty(x,y) = \abs{x} + \abs{y}$, this cross section is exactly $\abs{x}+1/2$. For $\phi_2(x,y) = \sqrt{x^2 + y^2}$, the cross section is a smooth curve, which cannot be captured perfectly by our discretization, but is better approximated than the cross section of $\phi_1(x,y) = \max\{\abs x, \abs y\}$, which has two kinks. The accuracy of the method depends on how well these cross sections can be approximated, since any error in these approximations will propagate to other regions.

With this in mind, we note that for $\phi_1(x,y) = \max\{\abs x , \abs y\}$, while the cross sections in the horizontal and vertical direction have these two kinks, the cross sections in the diagonal directions $y = x_0 \pm x$ will look like absolute value functions. If we used first-order approximations to $\nabla \phi_1$ along these diagonals, we would perfectly capture these cross sections, and thus reconstruct the solution exactly. This suggests that we should rotate the grid and consider alternative approximations to $\nabla \phi_1$. \Cref{sec:rotgrid} develops this idea. \\

\subsection{Increasing Accuracy with WENO Approximations} \label{sec:WENO}

We noted earlier that in some cases, one can increase the order of accuracy using (Weighted) Essentially Non-Oscillatory (WENO) schemes. The philosophy of ENO and WENO schemes---pioneered by Osher, Shu and Jiang, among others \cite{JiangWENO,OsherShu,Shu}---is to use multiple higher order approximations of $\phi_x$ and $\phi_y$, and deftly combine the approximations so as to minimize oscillations in the numerical solution near kinks. These methods were originally developed for time-dependent Hamilton-Jacobi equation, but have since been adapted to fast sweeping methods. We demonstrate the application of the third-order WENO approximations to our method, following the work of Zhang et al. \cite{ZhangWENOFSM}. One could use higher order WENO approximations if desired. 

The third-order WENO approximations to $\phi_x$ are given by \begin{equation}\label{eq:wenoapprox} \begin{split}
(\phi_x)_{ij}^{+} &= (1-w_x^+)\left( \frac{\phi_{i+1,j} - \phi_{i-1,j}}{2\Delta x}\right) + w_{x}^+ \left( \frac{-\phi_{i+2,j} + 4\phi_{i+1,j} - 3\phi_{ij}}{2\Delta x} \right),\\
(\phi_x)_{ij}^{-} &= (1-w_x^-)\left( \frac{\phi_{i+1,j} - \phi_{i-1,j}}{2\Delta x}\right) + w_{x}^- \left( \frac{\phi_{i-2,j} - 4\phi_{i-1,j} + 3\phi_{ij}}{2\Delta x} \right),
\end{split}
\end{equation} where the weights $w_x^+$ and $w_x^-$ are given by \begin{equation}\label{eq:weights}
\begin{split}
w^+_x &= \frac{1}{1+2(r_x^+)^2}, \,\,\,\,\,\,\,\,\,\,\,\,\,\, r_x^+ = \frac{\eps + (\phi_{i+2,j} - 2 \phi_{i+1,j} + \phi_{ij})^2}{ \eps + (\phi_{i+1,j} - 2 \phi_{ij} + \phi_{i-1,j})^2}, \\w^-_x &= \frac{1}{1+2(r_x^-)^2}, \,\,\,\,\,\,\,\,\,\,\,\,\,\, r_x^- = \frac{\eps + (\phi_{i-2,j} - 2 \phi_{i-1,j} + \phi_{ij})^2}{ \eps + (\phi_{i-1,j} - 2 \phi_{ij} + \phi_{i+1,j})^2}.
\end{split}
\end{equation} Here $\eps$ is some small number which we fix at $10^{-6}$. We define $(\phi_y)_{ij}^{+}$ and $(\phi_y)_{ij}^{-}$ analogously. 

Notice that each of the divided differences in \eqref{eq:wenoapprox} is a second-order approximation to $\phi_x$. The weighted averages---which favor the less oscillatory approximations---ensure that $(\phi_x)_{ij}^{+}$ and $(\phi_x)_{ij}^{-}$ are third-order approximations to $\phi_x$ in regions where $\phi$ is smooth. For a derivation and discussion of these formulas, see \cite{Shu} and the references therein. 

The question then becomes: how to include these approximations in a fast sweeping scheme? If we simply replace the finite difference approximations in \eqref{eq:upwindApprox} with $(\phi_x)^+_{ij}$ or $(\phi_x)^{-}_{ij}$ as appropriate, then we will not be able to isolate $\phi_{ij}$ and arrive at a simple update rule of the form \eqref{eq:updateRule}. The idea presented by Zhang et al.\cite{ZhangWENOFSM} is to start from the update rule itself. Note that the update rule \eqref{eq:updateRule} gives $\phi_{ij}$ as a function of $\phi_{i\pm1,j}$ and $\phi_{i,j\pm1}$. A finite difference approximation exploits the formal relationship $\phi(x\pm\Delta x,y) \approx \phi(x,y) \pm \Delta x \phi_x(x,y)$. Thus to arrive at a higher order approximation of the form \eqref{eq:updateRule}, we can  replace $\phi_{i+1,j}$ with $\phi_{ij} + \Delta x (\phi_{x})^+_{ij}$, and replace $\phi_{i-1,j}$ with $\phi_{ij} - \Delta x(\phi_x)^-_{ij}$, and similarly for $\phi_{i,j\pm 1}$. Doing so results in the update rule \begin{equation}\label{eq:wenoUpdateRule}
\phi_{ij}^*(a) = \frac{r_{ij} + \frac{\abs{f_{1,ij}(a)}}{\Delta x}\left(\phi_{ij} + \Delta x \xi_{1,ij}(a) (\phi_{x})^{\xi_{1,ij}(a)}_{ij}\right) + \frac{\abs{f_{2,ij}(a)}}{\Delta y}\left(\phi_{ij} + \Delta y \xi_{2,ij}(a) (\phi_{y})^{\xi_{2,ij}(a)}_{ij}\right)}{ \frac{\abs{f_{1,ij}(a)}}{\Delta x} +  \frac{\abs{f_{2,ij}(a)}}{\Delta y}}.
\end{equation} Using this update rule in \cref{alg:sweepingScheme} yields a higher order approximation of \eqref{eq:ourEquation2}. 

Formally, the approximation is third-order accurate when the solution $\phi(x,y)$ is smooth. In practice, the convergence can be corrupted by non-smoothness of the solution, and by the non-monotone nature of higher order approximations, which affects the numerical causality. Because of this last concern, when using the WENO approximations, it is crucial to seed the Gauss-Seidel iteration with a good initial guess $\phi^0_{ij}$, rather than simply setting $\phi_{ij}^0 = g_{ij}$ near the prescribed boundaries, and $\phi^0_{ij} = +\infty$ elsewhere. If one uses this crude initialization, it is easily checked in simple examples that \eqref{eq:wenoUpdateRule} will not correctly propagate information from the boundaries. We suggest first running the basic scheme with the ordinary first-order approximations, and using the resulting solution to initialize the iteration that uses the WENO approximations. 

We have carried out the implementation for two example problems. Both are of the form \begin{equation} 
r(x,y) = \|\nabla \phi(x,y)\|_2, \,\,\,\,\,\,\,\,\,\,\,\,\,\, \phi(0,0) = 0.
\end{equation} In the first, we take $r(x,y) = 1$ so that it is the same 2-norm Eikonal equation as above, and the solution is given by $\phi(x,y) = \sqrt{x^2 + y^2}$, which has a kink at the origin. In the second, we take $r(x,y) = \sqrt{x^2 + y^2}$, in which case the exact solution is $\phi(x,y) = (x^2 + y^2)/2$ which is smooth throughout the domain. The results are summarized in \cref{tab:weno1} and \cref{tab:weno2}. In this case we report both the $L^\infty$ and $L^1$ errors. In some cases, the $L^1$ error is more appropriate for evaluating the performance of WENO schemes, since the most significant errors can propagate along very small sets, whereas error remains small in the majority of the domain \cite{ZhangWENOFSM}. In \cref{tab:weno1}, we see that for the Eikonal equation $\| \nabla \phi \|_2 = 1$, the non-smoothness of the solutions corrupts the effects of the WENO approximations, and while the errors are smaller and convergence rate is improved, we do not nearly have third-order convergence. By contrast, in \cref{tab:weno2} when the solution remains smooth, we do see a greatly improved rate of convergence which is near third-order as the grid refines.

\begin{table}[t!]
\centering
\begin{subtable}[h]{\columnwidth}
\centering
\begin{tabular}{r r r r r}
$I,J$ & $L^\infty$ Err. & $L^\infty$ Conv. & $L^1$ Err. & $L^1$ Conv.\\ [0.5ex]
\hline \hline
50 & 4.3754e-02 &	--- & 9.7606e-02 	& ---\\
100 & 2.6310e-02 &	 0.7338 	& 5.9553e-02 & 0.7128   \\
200 & 1.5464e-02 &	 0.7666 	& 3.5451e-02 & 0.7484 \\
400 & 8.9201e-03 &	 0.7938 	& 2.0691e-02 & 0.7768\\ [0.5ex] \hline
\end{tabular}
\caption{Convergence table with first-order approximations to $\nabla \phi$.}
\end{subtable} 

\bigskip

\begin{subtable}[h]{\columnwidth}
\centering
\begin{tabular}{r r r r r}
$I,J$ & $L^\infty$ Err. & $L^\infty$ Conv. & $L^1$ Err. & $L^1$ Conv.\\ [0.5ex]
\hline \hline
50 & 9.0508e-03& 	    --- 	 & 2.0426e-02 	 &   ---\\
100  &4.4930e-03 &	     1.0104& 	 8.7373e-03 &	     1.2252 \\
200  & 2.2253e-03 &	     1.0137& 	 3.8868e-03 &	     1.1686 \\
400	& 1.0668e-03 &	     1.0607& 	 1.9013e-03 &	     1.0316 \\ [0.5ex] \hline
\end{tabular}
\caption{Convergence table with third-order WENO approximations to $\nabla \phi$.}
\end{subtable}
\caption{Error in solution of $\|\nabla \phi\|_2 = 1$ when using the Basic Method with (a) first-order approximations or (b) third-order WENO approximations. }
\label{tab:weno1}
\end{table}

\begin{table}[t!]
\centering
\begin{subtable}[h]{\columnwidth}
\centering
\begin{tabular}{r r r r r}
$I,J$ & $L^\infty$ Err. & $L^\infty$ Conv. & $L^1$ Err. & $L^1$ Conv.\\ [0.5ex]
\hline \hline
50 & 4.0010e-02 &	 --- &	 8.0016e-02 &	 ---\\
100 &	 2.0009e-02 &	 0.9997 &	 4.0014e-02 &	 0.9998  \\
200 &	 1.0010e-02 &	 0.9992 &	 2.0014e-02 &	 0.9995 \\
400 &	 5.0103e-03 &	 0.9985 &	 1.0014e-02 &	 0.9990\\ [0.5ex] \hline
\end{tabular}
\caption{Convergence table with first-order approximations to $\nabla \phi$.}
\end{subtable} 

\bigskip

\begin{subtable}[h]{\columnwidth}
\centering
\begin{tabular}{r r r r r}
$I,J$ & $L^\infty$ Err. & $L^\infty$ Conv. & $L^1$ Err. & $L^1$ Conv.\\ [0.5ex]
\hline \hline
50 & 	 2.3922e-03 &	     --- &	 5.4938e-03 & 	     --- \\
100  & 1.1609e-03 &	     1.0431 &	 2.3126e-03 & 	     1.2483  \\
200  & 	 1.5113e-04 &	     2.9413 &	 3.7584e-04 & 	     2.6213  \\
400	& 	 3.9126e-05 &	     1.9496 &	 6.0658e-05 & 	     2.6314  \\ [0.5ex] \hline
\end{tabular}
\caption{Convergence table with third-order WENO approximations to $\nabla \phi$.}
\end{subtable}
\caption{Error in solution of $\|\nabla \phi\|_2 = \sqrt{x^2 + y^2}$ when using the Basic Method with (a) first-order approximations or (b) third-order WENO approximations.}
\label{tab:weno2}
\end{table}

\section{A Rotating-Grid Fast Sweeping Scheme} \label{sec:rotgrid} In this section, we would like to append the basic algorithm with additional approximations to the gradient $\nabla \phi$ in directions that are not vertical and horizontal (with respect to the rectangular domain). In doing so, we can increase accuracy while maintaining a monotone scheme, since we do not use higher order approximations to the derivatives. 

In order to accomplish this, we must first recast equation \eqref{eq:ourEquation2d} in new coordinates $(\overline x, \overline y)$, rotated versions of the standard Cartesian coordinates. Again, we describe this procedure in two dimensions. Here the extension to higher dimensions is not as straightforward but can still be accomplished in a somewhat principled, if tedious, manner. We discuss the three-dimensional implementation in \cref{sec:append}. 

\begin{figure}[b!]
\centering
\includegraphics[width = 0.8\textwidth]{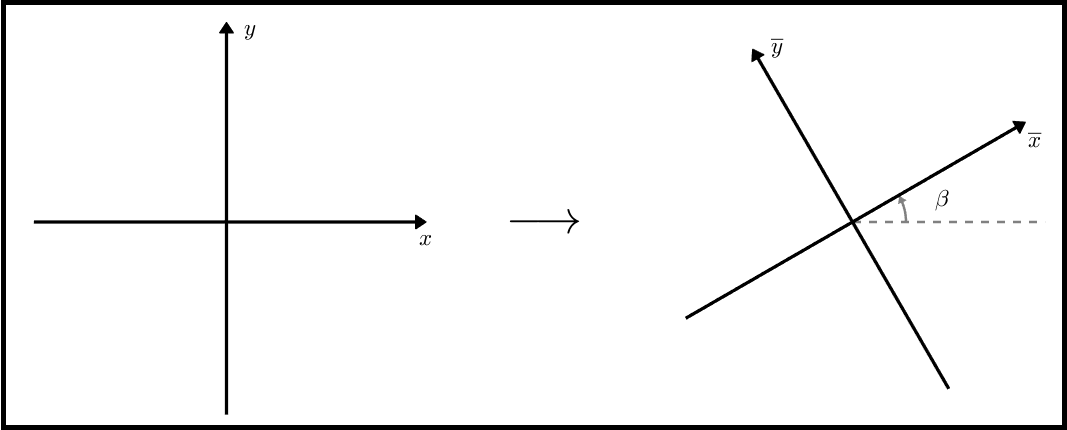}
\caption{Cartesian coordinates rotated by $\beta \in (0,\pi/2)$ in the counterclockwise direction.}
\label{fig:coordRotate}
\end{figure}

Suppose that $(\overline x, \overline y)$ are the typical Cartesian coordinates, rotated counterclockwise by an angle $\beta \in (0,\pi/2)$, as pictured in \cref{fig:coordRotate}. Note that it is sufficient to consider this range of angles; rotations by larger angles results in the same transformation up to renaming coordinates and flipping positive and negative directions. One easily verifies the relationship \begin{equation}\label{eq:newCoords} \begin{pmatrix} \overline x \\ \overline y \end{pmatrix} = \begin{pmatrix}\hphantom{-}\cos(\beta) & \sin(\beta) \\ -\sin(\beta) & \cos(\beta) \end{pmatrix} \begin{pmatrix}  x\\  y \end{pmatrix} \,\,\,\,\, \longleftrightarrow \,\,\,\,\,
\begin{pmatrix} x \\ y \end{pmatrix} = \begin{pmatrix}\cos(\beta) & -\sin(\beta) \\ \sin(\beta) & \hphantom{-}\cos(\beta) \end{pmatrix} \begin{pmatrix} \overline x\\ \overline y \end{pmatrix}.
\end{equation}
Thus the derivatives in the $(x,y)$ directions can be expressed \begin{equation}\label{eq:transformDeriv}
 \begin{split}
\phi_x &= \frac{\partial \overline x}{\partial x}~\phi_{\overline x}+ \frac{\partial \overline y}{\partial x}~\phi_{\overline y} = \cos(\beta) \phi_{\overline x} - \sin(\beta) \phi_{\overline y}, \\
\phi_y &= \frac{\partial \overline x}{\partial y}~\phi_{\overline x}+ \frac{\partial \overline y}{\partial y}~\phi_{\overline y} = \sin(\beta) \phi_{\overline x}+  \cos(\beta) \phi_{\overline y}.
\end{split} \end{equation} Inserting these representations into \eqref{eq:ourEquation2d} yields \begin{equation} \label{eq:ourEquation2dRot1} \begin{split}
-r(\overline x, \overline y) = \inf_{a\in A} \Big\{[\cos(\beta)&f_1(\overline x,\overline y,a) + \sin(\beta)f_2(\overline x,\overline y,a)] \phi_{\overline x}(\overline x,\overline y)\\ &+ [\cos(\beta)f_2(\overline x,\overline y,a) - \sin(\beta)f_1(\overline x,\overline y,a)] \phi_{\overline x}(\overline x,\overline y)  \Big\}. \end{split}
\end{equation} Defining \begin{equation}\label{eq:rotatedFuncs}\begin{split}
\overline f_1(\overline x,\overline y, a) &= \cos(\beta)f_1(\overline x,\overline y,a) + \sin(\beta)f_2(\overline x,\overline y,a), \\
\overline f_2(\overline x, \overline y, a) &= \cos(\beta)f_2(\overline x,\overline y,a) - \sin(\beta)f_1(\overline x,\overline y,a),
\end{split}\end{equation} we arrive at \begin{equation} \label{eq:ourEquation2dRot}
-r(\overline x, \overline y) = \inf_{a\in A} \big\{ \overline f_1(\overline x,\overline y,a) \phi_{\overline x}(\overline x,\overline y) + \overline f_2(\overline x,\overline y,a)\phi_{\overline y}(\overline x,\overline y) \big\}.
\end{equation} The idea is now to write the upwind finite difference approximations in the directions of $(\overline x, \overline y)$. Doing so shows that  \begin{equation} \label{eq:updateRuleRotated} 
\phi(\overline x,\overline y) =  \frac{r(\overline x, \overline y) + \frac{\abs{\overline f_1 (\overline x, \overline y,a)}}{\Delta \overline x} \phi(\overline x + \overline\xi_1 \Delta \overline x, \overline y) +  \frac{\abs{\overline f_2 (\overline x, \overline y,a)}} {\Delta \overline y} \phi(\overline x, \overline y + \overline\xi_2 \Delta \overline y)}{\frac{\abs{\overline f_1 (\overline x, \overline y,a)}}{\Delta \overline x} + \frac{\abs{\overline f_2 (\overline x, \overline y,a)}}{\Delta \overline y}}
\end{equation} is a first-order, upwind approximation to \eqref{eq:ourEquation2dRot} at the point $(\overline x, \overline y)$ when $a$ is the correct control value, and $\overline \xi_\ell = \text{sign}(\overline f_\ell (\overline x, \overline y,a))$. Thus one could add this approximation into the sweeping scheme and use the update rule \begin{equation}\label{eq:newUpdate}
\phi^n_{ij} = \min\Big\{\phi^{n-1}_{ij}, \min_{a\in A} \phi^*_{ij}(a), \min_{a\in A} \overline \phi^*_{ij}(a)\Big\},
\end{equation} where $\overline \phi^*_{ij}(a)$ is computed from \eqref{eq:updateRuleRotated}. However, this raises the question of how to evaluate \eqref{eq:updateRuleRotated} on the grid, since for example, $(\overline x \pm \Delta \overline x, \overline y)$ may not be grid nodes. 

Rotated finite differences are extensively used in computational wave mechanics. So-called rotated-staggered-grid methods were introduced by Saenger et al. \cite{RSG} and are still being developed and improved today \cite{RSG3,RSG2,RSG4,RSG5,RSG1}. The philosophy of these methods is the same: using finite differences in multiple orientations will more accurately capture the upwind direction. Their strategy is to define a new grid corresponding to the points $(\overline x, \overline y)$ and keep track of solution values $\phi_{ij}$ and $\overline \phi_{ij}$ separately, while using both sets of values to approximate the derivatives on both grids. To this author's knowledge, the idea of fixing a square grid and computing approximations to $\nabla \phi$ in different directions has not been widely used in the context of fast sweeping methods. Takei et al. \cite{TakeiTsai1} suggest using approximations along different directions. However, in their case, the upwind direction is fixed (in analogy to our setup, they have $f_1, f_2$ independent of $a$) which simplifies the matter.

We would like to maintain a single grid $(x_i,y_j)$. To do so, one could interpolate values of $\phi_{ij}$ to off grid values, and compute the upwind approximation in any direction $\beta$. This would be computationally expensive since, in order to maintain the Gauss-Seidel sweeping, this interpolation will need to be performed separately for every $(i,j)$ using the newest updated values. Alternatively, we can choose particular values of $\beta$ and $\Delta \overline x, \Delta \overline y$ such that the points $(\overline x \pm \Delta \overline x, \overline y)$, $(\overline x, \overline y \pm \Delta \overline y)$ fall on the grid. 

\begin{figure}[t!]
\centering
\includegraphics[width = 0.6\textwidth,trim=40 40 40 40,clip]{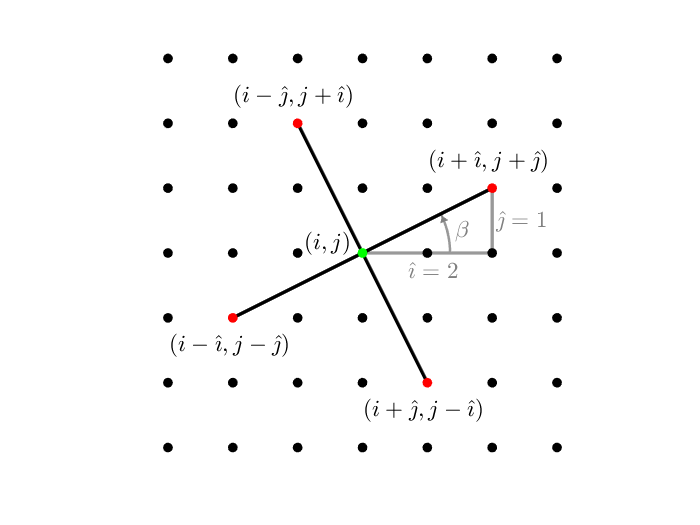}
\caption{Rotated stencil at $(i,j)$ using the rotation determined by $(\hat \imath,\hat\jmath) = (2,1)$.}
\label{fig:gridRotate}
\end{figure}

Explicitly, rather than choosing $\beta$ and the rotated grid parameters $(\Delta \overline x, \Delta \overline y)$, we choose natural numbers $(\hat\imath,\hat\jmath)$, and define $\beta = \arctan(\hat\jmath/\hat\imath)$. We then let this $\beta$ determine the grid rotation, so that the positive $\hat x$-direction is parallel with the vector $(\hat \imath, \hat \jmath)$. This is pictured in \cref{fig:gridRotate}. Here we have used $(\hat\imath,\hat\jmath) = (2,1)$. As pictured, the nodes used to approximate $\phi_{\overline x}$ at $(i,j)$ will be $\{(i,j), (i+2,j+1)\}$ for the forward approximation, and $\{(i-2,j-1),(i,j)\}$ for the backward approximation. Similarly, the nodes used to approximate $\phi_{\overline y}$ at $(i,j)$ will be $\{(i,j), (i-1,j+2)\}$ for the forward approximation and $\{(i+1,j-2),(i,j)\}$ for the backward approximation.

We note that as described, this will only work on a square grid $(\Delta x = \Delta y)$. The extension to a non-square grid is a bit more complicated. In that case, there would be two rotation angles that rotate the $x$-axis and $y$-axis differently, and thus the resulting coordinate system would no longer be orthogonal. For the remainder of this document, we will assume that $\Delta x = \Delta y$ so that the rotation method works as described.

\begin{figure}[t!]
\centering
\begin{subfigure}{0.6\textwidth}
\centering 
\includegraphics[width=\textwidth,trim = 40 40 40 40,clip]{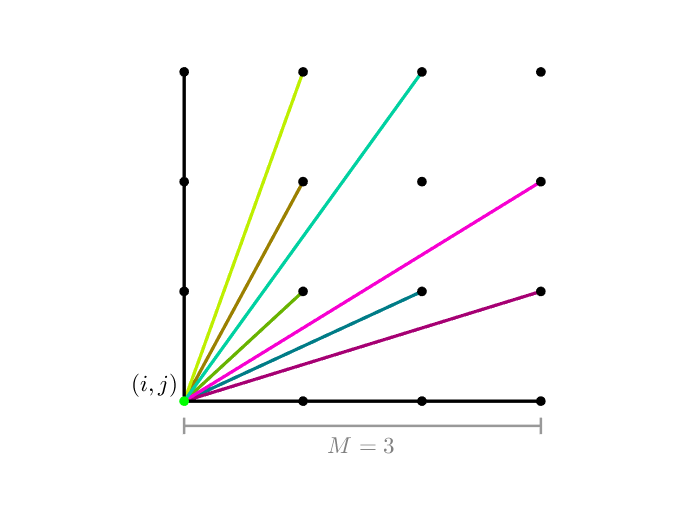}
\caption{The possible rotation angles if $1 \le \hat\imath,\hat\jmath \le 3$.}
\end{subfigure}~
\begin{subfigure}{0.3\textwidth}
\centering
\begin{tabular}{r r}
$M$ & $C(M)$ \\ [0.5ex]
\hline \hline 
  1& 	 1\\
  2& 	 3\\
 3& 	 7\\
 4& 	 11\\
 5& 	 19\\
6& 	23\\ 
7& 35\\
8 & 43\\
9& 55\\
10& 63\\[0.5ex]
\hline 
\end{tabular}
\caption{Number of possible \hphantom{hi} angles if $1 \le \hat\imath,\hat\jmath \le M$.}
\end{subfigure}
\caption{If we restrict $1 \le \hat\imath,\hat\jmath \le M$ there will be some finite number $C(M)$ of distinct rotation angles $\beta = \arctan(\hat\jmath/\hat\imath)$, each represented by a colored line.}
\label{fig:differentAngleRotations}
\end{figure}

With these parameters $(\hat\imath,\hat\jmath)$ determining the rotation, we define the new grid discretization parameter $\Delta s = \sqrt{(\hat\imath \Delta x)^2 + (\hat\jmath \Delta y)^2}$. Note that this $\Delta s$ will take the place of $\Delta \overline x, \Delta \overline y$ in the case of a square grid. Thus we can translate equation \eqref{eq:updateRuleRotated} onto the grid: \begin{equation}\label{eq:updateRuleRotatedGrid}
\overline \phi^*_{ij}(a) = \frac{r_{ij}\Delta s + \abs{\overline f_{1,ij}(a)} \phi_{i+\overline \xi_{1,ij}(a)\hat\imath, j+\overline\xi_{1,ij}(a)\hat\jmath} + \abs{\overline f_{2,ij}(a)} \phi_{i-\overline \xi_{2,ij}(a)\hat\jmath, j+\overline\xi_{2,ij}(a)\hat\imath}}{\abs{\overline f_{1,ij} (a)} + \abs{\overline f_{2,ij}(a)}},
\end{equation} which, one sees, is exactly analogous to \eqref{eq:updateRule}, except that the coordinates are rotated and the grid parameters are equal. Inserting this approximation into \eqref{eq:newUpdate} provides a new update rule that can be used in \cref{alg:sweepingScheme}. Of course, it is not necessary to limit oneself to a single rotation $(\hat\imath,\hat\jmath)$. To further improve the scheme, one can choose as many pairs $(\hat \imath, \hat \jmath)$ as desired, compute the rotated derivative approximations in each of these directions, and take the minimum over all such approximations. Since the stencil at each grid node will be larger, the scheme will require a larger layer of ghost nodes padding the computational boundary; otherwise, \cref{alg:sweepingScheme} will operate in the exact same fashion, but with extra approximations included in the update rule. In general, if one imposes $1 \le \hat\imath, \hat\jmath \le M$, one should buffer the computational domain with $M$ layers of grid nodes, and there will be some finite number $C(M)$
of distinct angles $\beta$ created by different pairs $(\hat\imath,\hat\jmath).$\footnote{In fact, one has $C(M) = 2\Big( \sum^M_{m=1} \varphi(m)\Big) -1$ where $\varphi$ is the Euler totient function, as detailed in the Online Encyclopedia of Integer Sequences: \url{http://oeis.org/A018805}}
~This is pictured in \cref{fig:differentAngleRotations}, where each colored line represents a distinct rotation angle $\beta$ when $M=3$. Fixing $M$, we propose two strategies for choosing different rotation angles: first, one could simply use every possible rotation angle. This may be computationally expensive since, for example, when $M = 5$, there are $C(M) = 19$ angles to consider. Accordingly, our second strategy will be to choose some fixed size subcollection at random. This will not be able to guarantee the same level of accuracy, but will be significantly cheaper computationally.  It may also be better than choosing a fixed subcollection of angles since, in application, one may not be able to intuit the ``principal" directions that need to be captured as we can for the Eikonal equations. Another possibility would be to change the rotation angle $\beta$ for each grid point, perhaps accounting for the admissible control actions and possible upwind directions; this is essentially what is done by rotating the grid by $\beta = \pi /4$ for the 1-norm Eikonal equation below. To do so more generally, one would need to carefully analyze the particular update rule \eqref{eq:updateRuleRotatedGrid} for one's problem in order to determine a range of possible upwind directions. As presented, we fix the rotation angles $\beta$ before each iteration. 

Note that we will always use the ordinary forward and backward approximations in the $(x,y)$ directions, and include approximations in other directions as desired. This is to establish a baseline. In this manner, using derivative approximations in additional directions can only improve upon the accuracy of the basic method presented in \cref{alg:sweepingScheme}. 

It is natural to consider the optimal number of grid rotations---or similarly, the optimal width of a stencil---for a given problem. Unfortunately, it is difficult to address this point generally. In specific examples, the answer is simple. For example, in the 1-norm Eikonal equation, one can achieve an exact solution with a single grid rotation, as we demonstrate in the succeeding section, and thus additional rotations will offer no benefit. However, for the 2-norm Eikonal equation, each new rotation will serve to better capture the solution at certain points, since characteristics travel outward from the origin in every direction. For general steady-state HJB equations, one may not know the characteristic directions ahead of time, so while adding more rotations can do no worse than the basic scheme, the benefits may be marginal, and they come at the cost of increasing the computational burden. Accordingly, this point would need to be addressed on an \emph{ad hoc} basis, and depends both on the problem and on the user's desire to balance the possibility of large accuracy gains against the increased computation.

\subsection{Application of the Rotating-Grid Method to Eikonal Equations} \label{sec:eikonalRot}
\begin{figure}[t!]
\centering
\begin{subfigure}{0.45\textwidth}
\centering 
\includegraphics[width=\textwidth,trim=43 15 35 5,clip]{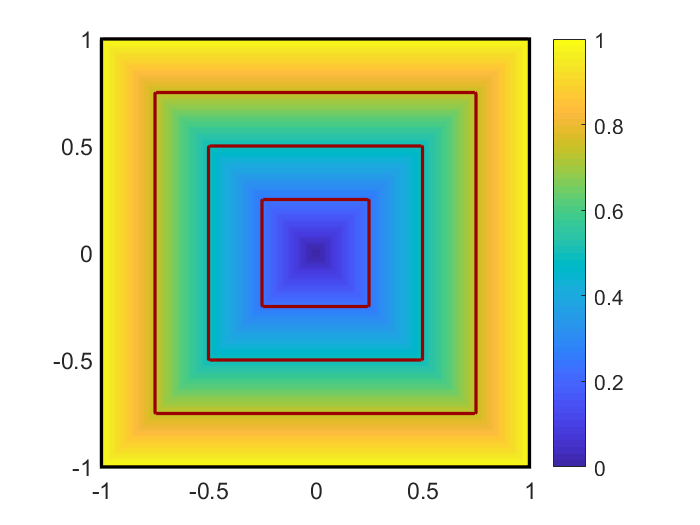}
\caption{Approximate solution}
\label{fig:OneNormSolnWithRot}
\end{subfigure}~
\begin{subfigure}{0.45\textwidth}
\centering 
\includegraphics[width=\textwidth,trim=43 15 35 5,clip]{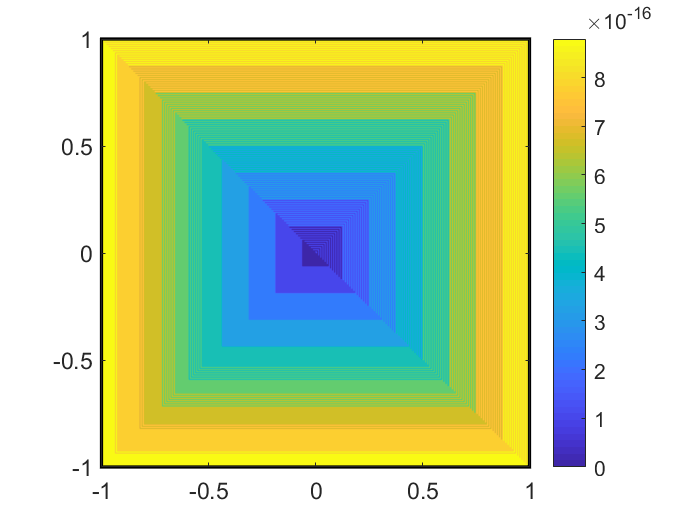}
\caption{Error in approximation.}
\label{fig:OneNormErrWithRot}
\end{subfigure}
\caption{Numerical solution of $\|\nabla \phi\|_1 = 1$ with additional approximations to $\nabla \phi$ in the direction of $\beta = \pi/4$. Compare with figures \ref{fig:OneNormSoln},~\ref{fig:OneNormErr}.}
\label{fig:OneNormWithRot}
\end{figure}

We apply the sweeping scheme with rotated derivative approximations to the Eikonal equation in the $p=1$ and $p=2$ norms. We remarked earlier that cross sections of the solution $\phi_1(x,y) = \max\{\abs{x},\abs y\}$ along the diagonal lines $y = x_0 \pm x$ could be captured exactly by our scheme if we use the rotation $\beta = \pi/4$, which is the same as $(\hat\imath,\hat\jmath) = (1,1)$. In this case, the rotated coefficients are $\overline f_1 = \frac 1 {\sqrt 2}(a_1+ a_2)$ and $\overline f_2 = \frac 1 {\sqrt 2} (a_2-a_1)$, where $a_1,a_2\in \{\pm 1\}$. Since one of these is zero, the update rule is \begin{equation} \label{eq:updateRule1normWithRotation} \begin{split}
\phi^{n,1}_{ij} = \min\bigg\{  \phi^{n-1,1}_{ij}, &\frac{1 + \frac{1}{\Delta x} \phi^{n,1}_{i+1,j} + \frac{1}{\Delta y} \phi^{n,1}_{i,j+1}}{\frac{1}{\Delta x} + \frac 1 {\Delta y}},  \frac{1 + \frac{1}{\Delta x} \phi^{n,1}_{i-1,j} + \frac{1}{\Delta y} \phi^{n,1}_{i,j+1}}{\frac{1}{\Delta x} + \frac 1 {\Delta y}},\\
& \frac{1 + \frac{1}{\Delta x} \phi^{n,1}_{i+1,j} + \frac{1}{\Delta y} \phi^{n,1}_{i,j-1}}{\frac{1}{\Delta x} + \frac 1 {\Delta y}}, \frac{1 + \frac{1}{\Delta x} \phi^{n,1}_{i-1,j} + \frac{1}{\Delta y} \phi^{n,1}_{i,j-1}}{\frac{1}{\Delta x} + \frac 1 {\Delta y}},\\
&\phi^{n,1}_{i+1,j+1} + \frac{\Delta s}{\sqrt 2},\,\, \phi^{n,1}_{i-1,j-1} + \frac{\Delta s}{\sqrt 2}\\
&\phi^{n,1}_{i-1,j+1} + \frac{\Delta s}{\sqrt 2},\,\,\phi^{n,1}_{i+1,j-1} + \frac{\Delta s}{\sqrt 2}\bigg\}.
\end{split} 
\end{equation} We use this update rule in \cref{alg:sweepingScheme} to solve $\|\nabla \phi\|_1=1$. The results are seen \cref{fig:OneNormWithRot}. We note that the level sets of the solution have sharp edges, as opposed to \cref{fig:OneNormSoln}, where they were rounded off. In this case, the error in the solution is on the order of machine-$\eps$.

\begin{figure}[t!]
\begin{subfigure}{0.32\textwidth}
\centering 
\includegraphics[width=\textwidth,trim=43 15 35 5,clip]{images/TwoNormSoln}
\caption{Approx. soln., no additional directions.}
\label{fig:TwoNormSolnComp}
\end{subfigure}\hfill
\begin{subfigure}{0.32\textwidth}
\centering 
\includegraphics[width=\textwidth,trim=43 15 35 5,clip]{images/TwoNormErr}
\caption{Error in approx., no additional directions.}
\label{fig:TwoNormErrComp}
\end{subfigure}\hfill
\begin{subfigure}{0.35\textwidth}
\centering
\begin{tabular}{r r r}
$I,J$ & $L^\infty$ Err. & Conv.  \\ [0.5ex]
\hline \hline 
  50& 	 4.3754e-02& 	 ---\\
 100& 	 2.6310e-02& 	 0.7338\\
 200& 	 1.5464e-02& 	 0.7666\\
 400& 	 8.9201e-03& 	 0.7938\\
 800& 	 5.0668e-03& 	 0.8160\\
1600& 	 2.8431e-03& 	 0.8336\\ [0.5ex]
\hline 
\end{tabular}
\caption{Conv. table, no additional directions.}
\label{fig:TwoNormConvComp}
\end{subfigure} \\
\begin{subfigure}{0.32\textwidth}
\centering 
\includegraphics[width=\textwidth,trim=43 15 35 5,clip]{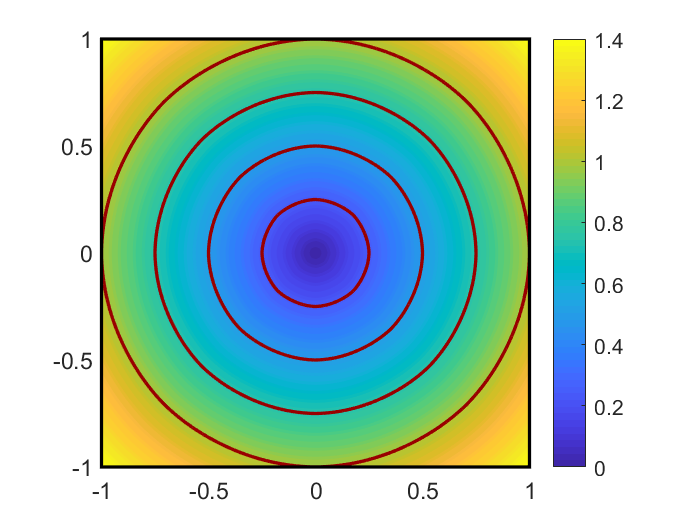}
\caption{Approx. soln., one additional direction.}
\label{fig:TwoNormSolnOneDirection}
\end{subfigure}\hfill
\begin{subfigure}{0.32\textwidth}
\centering 
\includegraphics[width=\textwidth,trim=43 15 35 5,clip]{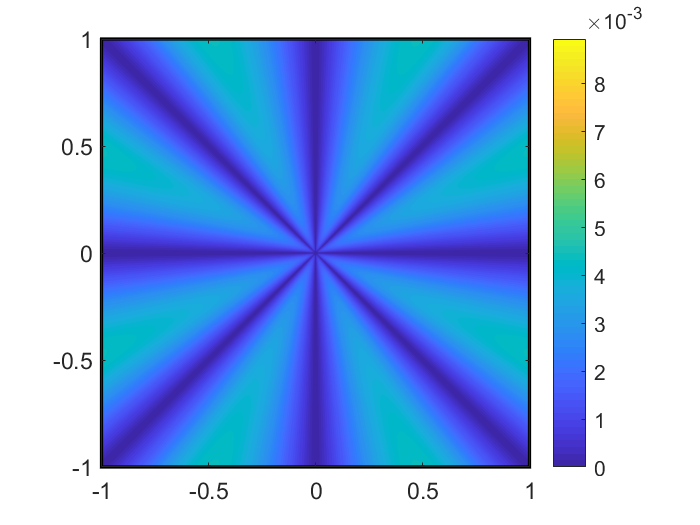}
\caption{Error in approx., one additional direction.}
\label{fig:TwoNormErrOneDirection}
\end{subfigure}\hfill
\begin{subfigure}{0.35\textwidth}
\centering
\begin{tabular}{r r r}
$I,J$ & $L^\infty$ Err. & Conv.  \\ [0.5ex]
\hline \hline 
  50& 	 1.7901e-02& 	 ---\\
 100& 	 1.1567e-02& 	 0.6300\\
 200& 	 7.2269e-03& 	 0.6789\\
 400& 	 4.3888e-03& 	 0.7192\\
 800& 	 2.6063e-03& 	 0.7518\\
1600& 	 1.5202e-03& 	 0.7777\\ [0.5ex]
\hline 
\end{tabular}
\caption{Conv. table, one additional direction.}
\label{fig:TwoNormConvOneDirection}
\end{subfigure} \\
\begin{subfigure}{0.32\textwidth}
\centering 
\includegraphics[width=\textwidth,trim=43 15 35 5,clip]{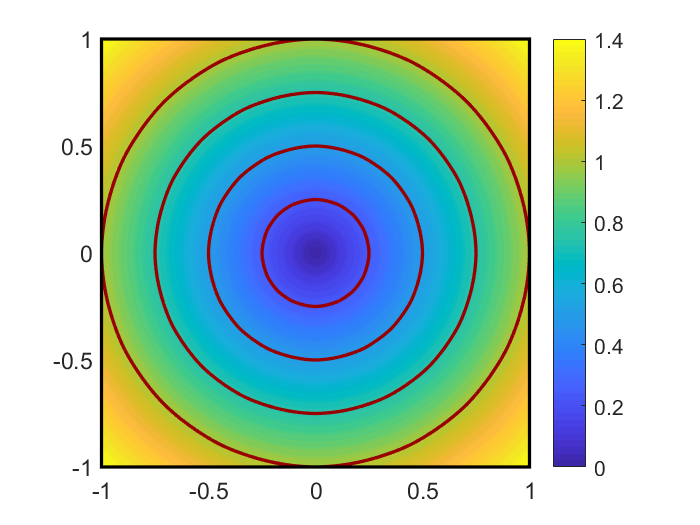}
\caption{Approx. soln., three additional directions.}
\label{fig:TwoNormSolnThreeDirections}
\end{subfigure}\hfill
\begin{subfigure}{0.32\textwidth}
\centering 
\includegraphics[width=\textwidth,trim=43 15 35 5,clip]{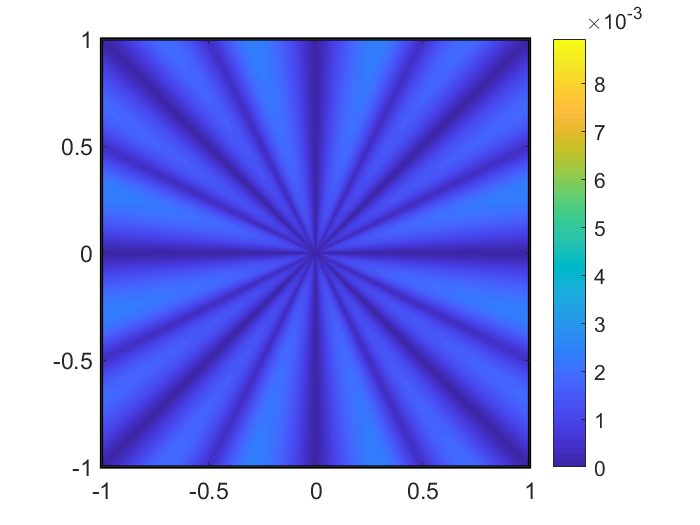}
\caption{Error in approx., three additional directions.}
\label{fig:TwoNormErrThreeDirections}
\end{subfigure}\hfill
\begin{subfigure}{0.35\textwidth}
\centering
\begin{tabular}{r r r}
$I,J$ & $L^\infty$ Err. & Conv.  \\ [0.5ex]
\hline \hline 
  50& 	 8.7787e-03& 	 ---\\
 100& 	 5.9351e-03& 	 0.5647\\
 200& 	 3.8508e-03& 	 0.6241\\
 400& 	 2.4134e-03& 	 0.6741\\
 800& 	 1.4720e-03& 	 0.7133\\
1600& 	 8.7876e-04& 	 0.7443\\[0.5ex]
\hline 
\end{tabular}
\caption{Conv. table, three additional directions.}
\label{fig:TwoNormConvThreeDirections}
\end{subfigure} 
\caption{Numerical solution of $\| \nabla \phi\|_2 = 1$ using our fast sweeping method with additional approximations to $\nabla \phi$ in different directions. Scale on error plots is fixed. Error is approximately zero in the directions of the derivative approximations.}
\label{fig:TwoNormComparison}
\end{figure}

Next we solve $\|\nabla \phi\|_2 = 1$. Here, in contrast with $\|\nabla \phi\|_1 = 1$ or $\|\nabla \phi\|_\infty = 1$, we will never be able to solve the equation exactly with finitely many grid rotations. The solution will be resolved exactly along any line through the origin if we consider the derivatives in the direction along that line. We saw this in \cref{fig:TwoNormErr}; the error is approximately zero along the $x$-axis and $y$-axis. We see it further in \cref{fig:TwoNormComparison}. In that figure, we first solve $\|\nabla \phi\|_2=1$ using the basic method (subfigures \ref{fig:TwoNormSolnComp},~\ref{fig:TwoNormErrComp},~\ref{fig:TwoNormConvComp}). We then compare this to results when using approximations to the derivatives in one additional direction (subfigures \ref{fig:TwoNormSolnOneDirection},~\ref{fig:TwoNormErrOneDirection},~\ref{fig:TwoNormConvOneDirection}), and three additional directions (subfigures \ref{fig:TwoNormSolnThreeDirections},~\ref{fig:TwoNormErrThreeDirections},~\ref{fig:TwoNormConvThreeDirections}). As expected, we see that for a fixed $I,J$, the error only decreases as we incorporate additional appoximations to $\nabla \phi$ in different directions. Interestingly, the order of convergence appears to slightly decrease when additional directions are included. However, we also note that when using three additional directions one only needs $51$ grid points in each direction to achieve the same approximation error as the basic method with $401$ points in each direction. 

\begin{figure}[t!]
\centering
\begin{subfigure}{0.48\textwidth}
\centering 
\includegraphics[width=\textwidth,trim=43 15 35 5,clip]{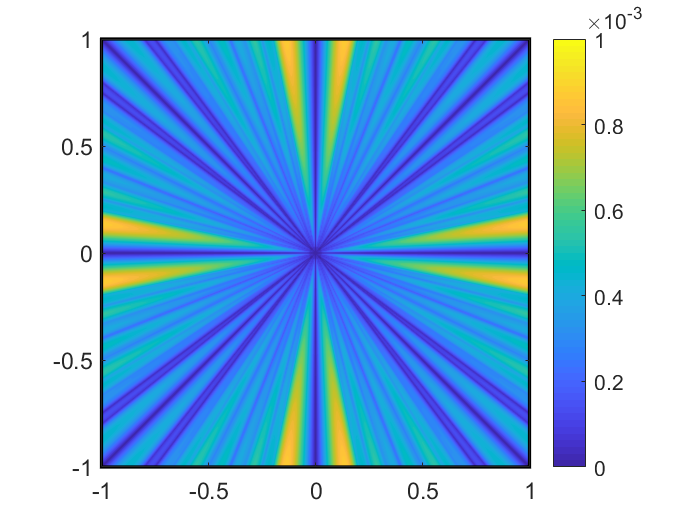}
\caption{Error in approximation when all 19 grid rotations are used in each iteration. Maximum error is $8.7914 \times 10^{-4}$. \\ \hphantom{hi}}
\label{fig:TwoNormErrStencil5}
\end{subfigure}\hfill
\begin{subfigure}{0.48\textwidth}
\centering 
\includegraphics[width=\textwidth,trim=43 15 35 5,clip]{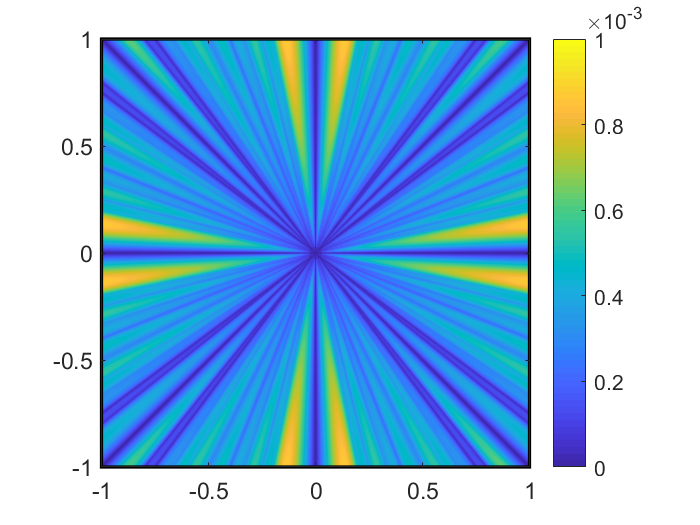}
\caption{Error in approximation when each iteration uses 2 grid rotations chosen randomly from the 19 possibilities. Maximum error is $8.7941 \times 10^{-4}$}
\label{fig:TwoNormErrStencil5Rand}
\end{subfigure}
\caption{Error in approximation using rotations $\beta = \arctan(\hat\jmath/\hat \imath)$ where $1\le \hat \imath,\hat\jmath \le 5$.}
\label{fig:TwoNormErrCompRand}
\end{figure}

Finally, we solve the same equation using a $401\times 401$ grid and all 19 grid rotations $\beta = \arctan (\hat\jmath/\hat\imath)$ corresponding to $1 \le \hat \imath, \hat \jmath \le 5$. In \cref{fig:TwoNormErrStencil5}, we see that when using all 19 rotations, we achieve an approximation error of $8.7914 \times 10^{-4}$. In this case, the algorithm required 12 iterations to terminate, and each iteration requires 20 times the computation as in the basic method (since there are 20 total approximations to $\nabla \phi$ being computed). In \cref{fig:TwoNormErrStencil5Rand}, we use the same 19 possible grid rotations, but for each iteration we choose only two rotations to use at random. We achieve similar approximation error: $8.7941 \times 10^{-4}$. The algorithm required 40 iterations to converge, but each iteration is  3 times as costly as in the basic method. Thus while there are roughly 3 times as many iterations, each iteration requires only 15\% of the computation, meaning one can achieve similar approximation error with roughly half the computation. It should be mentioned that these results have some randomness, but the numbers presented are quite typical.

We note that Darbon and Osher \cite{Darbon} solve similar Eikonal equations using a variational method based on the Hopf-Lax formula. Their method is applicable in high dimensions and can resolve the solution with essentially no error. However, the method only applies to Hamiltonians which are state-independent: $H = H(\nabla \phi)$. Fast sweeping methods are more general, but suffer from the curse of dimensionality. We have included Eikonal equations as an example because they are the prototypical steady-state Hamilton-Jacobi equations.\\

\subsection{Iteration Counts and Comparison with the Lax-Friedrichs Sweeping Scheme}

One final consideration when weighing the efficiency of a sweeping scheme is the iteration count necessary for the scheme to converge. Accordingly, we include a brief discussion regarding the iteration counts for the algorithm with different derivative approximations. We note again that one of the primary strengths of our algorithm is its ease of implementation. One other fast sweeping method which shares this ease of implementation is the Lax-Friedrichs (LF) sweeping scheme devised by Kao, Osher and Qian \cite{Kao2004}. In two dimensions, their scheme approximates the equation $H(x,y,\phi_x,\phi_y) = r(x,y)$ using the update rule \begin{equation}\label{eq:LFUpdateRule}
\phi^*_{ij} = \frac{r_{ij} - H\left(x_i,y_j, \frac{\phi_{i+1,j} - \phi_{i-1,j}}{2\Delta x},\frac{\phi_{i,j+1} - \phi_{i,j-1}}{2\Delta y} \right) +  \sigma_x \frac{\phi_{i+1,j} + \phi_{i-1,j}}{2\Delta x}  +  \sigma_y \frac{\phi_{i,j+1} + \phi_{i,j-1}}{2\Delta y}}{\frac{\sigma_x}{\Delta x} + \frac{\sigma_y}{\Delta y}}.
\end{equation} Intuitively, one arrives at this formula by using the centered difference approximations to $\phi_x$ and $\phi_y$, and adding artificial viscosity at strength $O(\Delta x, \Delta y)$. Here $\sigma_x$ and $\sigma_y$ are the artificial viscosity coefficients; they are bounds on $\partial H / \partial \phi_x$ and $\partial H / \partial \phi_y$ respectively. 

This method applies to general steady-state Hamilton-Jacobi equations, and is easily implemented regardless of how complicated the Hamiltonian may be. This is in contrast to other fast sweeping schemes, wherein the local update rule entails solving a nonlinear equation whose complexity depends on the Hamiltonian \cite{Qian1,Qian2,ZhangWENOFSM}. The tradeoff is that due to the diffusive nature of the LF numerical Hamiltonian, there is no causality condition being enforced, and consequently, a very large number of iterations are required for convergence.

We demonstrate this using the 2-norm Eikonal equation $\| \nabla \phi \|_2 = 1$ on $[-1,1] \times [-1,1]$. Note that because the characteristics are straight lines flowing out of the origin, our basic scheme, being fully upwind, converges in a single iteration. When we include additional approximations to the derivatives in rotated directions, this is no longer true. The scheme is still upwind, but there are multiple approximations to a given derivative which obey the causality condition, and alternate iterations may prefer different approximations, which means the algorithm requires more than one iteration to converge. The results are contained in \cref{tab:LFcomp}. As seen in the table, the LF sweeping scheme requires significantly more iterations in order to converge, and results in a larger $L^\infty$ error. As expected, the basic method converges in one iteration for any grid resolution. If we add derivative approximations in different directions, the algorithm no longer converges in one iteration, but empirically, we notice that when we add more approximations, fewer iterations are required. In all of these tests, the convergence criterion is $\max_{ij} \abs{\phi^n_{ij} - \phi^{n-1}_{ij}} < 10^{-8}$.

It should be noted that, while the LF scheme requires more iterations, each iteration is more efficient since there is no minimization problem or nonlinear inversion. The LF schemes also applies to more general problems. However, in cases where the minimization in our scheme is easily resolved, it is likely to outperform the LF scheme both in terms of efficiency and accuracy. We see this with the last example in \cref{sec:otherApps}. 

\begin{table}[h!] 
\centering
\begin{tabular}{r | r r | r r | r r | r r }
& \multicolumn{2}{c}{Lax-Friedrichs} & \multicolumn{2}{c}{Basic} & \multicolumn{2}{c}{Basic$+1$} & \multicolumn{2}{c}{Basic$+3$}\\
\hline
$I,J$  & Iter. & $L^\infty$ Err. & Iter. & $L^\infty$ Err. & Iter. & $L^\infty$ Err. & Iter. & $L^\infty$ Err.\\ [0.5ex] \hline \hline
50 & 34&	1.0958e-01	&1	&4.3754e-02 & 5	&1.7901e-02 &	5 &	8.7787e-03\\
100	&	43 &	6.1799e-02 &	1 &	2.6310e-02 & 8 &	1.1567e-02 &	7 &	5.9351e-03 \\
200	&	59	&3.4387e-02	&1	&1.5464e-02 & 14 & 	7.2252e-03 &	10	& 3.8508e-03\\
400	&	91 &	1.8932e-02 &	1 &	8.9201e-03 & 24 &	4.3888e-03 &	18	& 2.4134e-03\\ \hline
\end{tabular} 
\caption{The iteration counts for different versions of our algorithm and for the Lax-Friedrichs sweeping scheme when solving $\| \nabla \phi\|_2 =1$. Here Basic$+1$ designates the basic method appended with derivative approximations in one additional direction; Basic$+3$ designates the basic method appended with derivative approximations in three additional directions.}
\label{tab:LFcomp}
\end{table}

\section{Other Applications}\label{sec:otherApps}

Lastly, we present two applications of our method to problems arising in engineering. First we consider the visibility problem. Here one could imagine placing cameras at fixed points in a domain. The cameras have omnidirectional view, but the view is occluded by obstacles. The problem is to find the region that is visible to the cameras. 

This problem was first formulated using partial differential equations and the level set method by Tsai et al. \cite{Tsai}. However, that formulation involves a nonlocal equation. More recently, Oberman and Salvador were able to recast the problem in terms of a simple, local equation \cite{ObermanVisibility}. Specifically, supposing that $g: \R^d \to \R$ is the signed distance function to the obstacles (positive inside the obstacles) and $x^*\in \R^d$ is the vantage point, the visibility function $\phi: \R^d \to \R$ satisfies \begin{equation}\label{eq:visibilityPDE}
0 = \min\{\phi(x) - g(x),~ \langle x-x^*, \nabla \phi(x)\rangle\}
\end{equation} with the boundary condition $\phi(x^*) = g(x^*)$. The visibility set is then given by $\{\phi \le 0\}$. To include multiple vantage points, one solves \eqref{eq:visibilityPDE} individually for each point, and combines the solution via minima and maxima to account for different scenarios (for example, the minimum of all such solutions will provide the set of points visible from at least one vantage point, while the maximum of all such solutions provides the set of points that are visible from all vantage points simultaneously).

\begin{figure}[b!]
\centering
\begin{subfigure}{0.4\textwidth}
\centering 
\includegraphics[width=\textwidth,trim=43 15 35 15,clip]{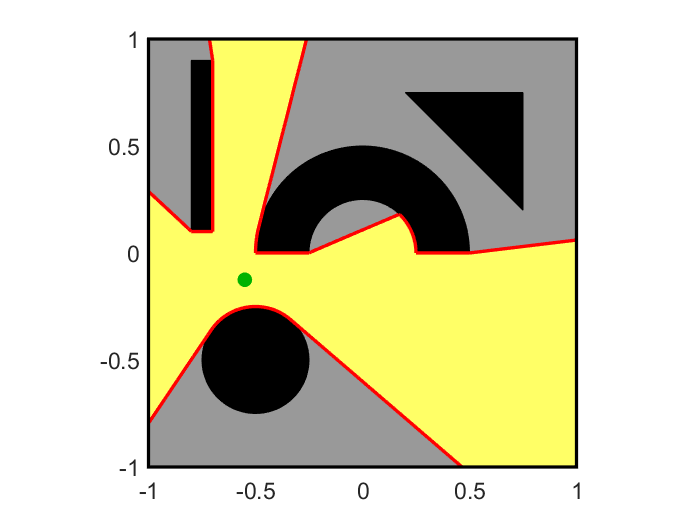}
\caption{One vantage point.}
\label{fig:Visibility1}
\end{subfigure}~
\begin{subfigure}{0.4\textwidth}
\centering 
\includegraphics[width=\textwidth,trim=43 15 35 15,clip]{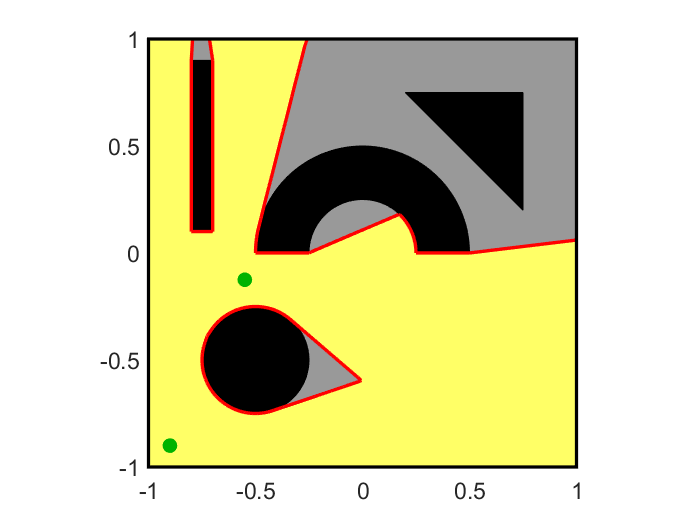}
\caption{Two vantage points.}
\label{fig:Visibility2}
\end{subfigure}\\
\begin{subfigure}{0.4\textwidth}
\centering 
\includegraphics[width=\textwidth,trim=43 15 35 15,clip]{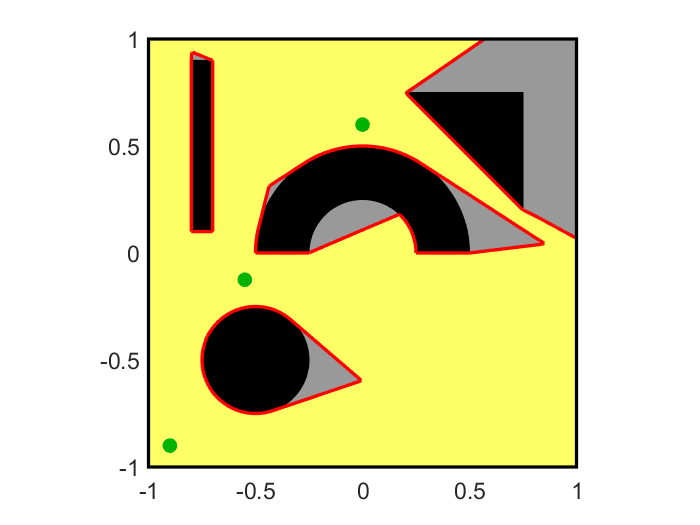}
\caption{Three vantage points.}
\label{fig:Visibility3}
\end{subfigure}~
\begin{subfigure}{0.4\textwidth}
\centering 
\includegraphics[width=\textwidth,trim=43 15 35 15,clip]{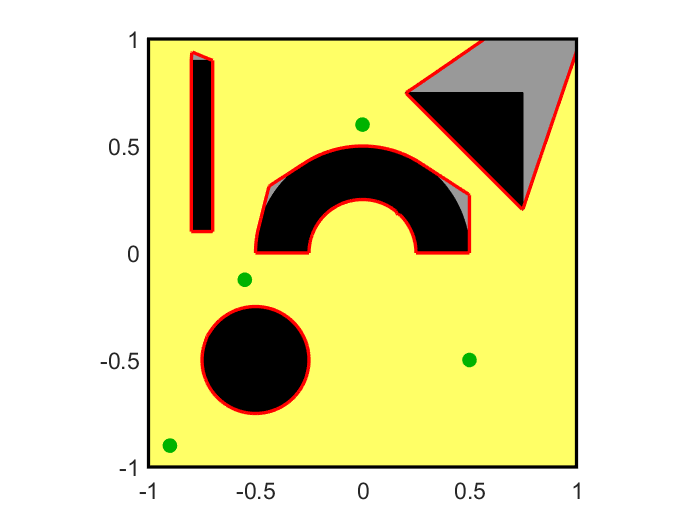}
\caption{Four vantage points.}
\label{fig:Visibility4}
\end{subfigure}
\caption{Computing the visibility set using \eqref{eq:visibilityPDE}. The green dots represent the vantage points. The black shapes are obstacles. The yellow set is comprised of points visible from at least one vantage point. The grey set is the unobserved set.}
\label{fig:VisibilityExample}
\end{figure}

Note that while equation \eqref{eq:visibilityPDE} does not directly follow from an optimal control problem, it does fit into our framework. If one sets $\phi^0_{ij} = g_{ij}$ for the nodes closest to the vantage point $(x^*,y^*)$ and $\phi^0_{ij} = -\infty$ at other nodes, one can use the update rule \begin{equation}\label{eq:updateRuleVisibility}
\phi^*_{ij} = \frac{\frac{\abs{x_i-x^*}}{\Delta x} \phi_{i-\text{sign}(x_i-x^*),j} + \frac{\abs{y_j-y^*}}{\Delta y} \phi_{i,j-\text{sign}(y_i-y^*)}}{\frac{\abs{x_i-x^*}}{\Delta x} + \frac{\abs{y_j-y^*}}{\Delta y}},
\end{equation} and iterate $\phi^n_{ij} = \max\{\phi^{n-1}_{ij},g_{ij},\phi^*_{ij}\}$. [Note that the upwind direction is reversed, which explains the slight deviations between these formulas and those above.] One can then use additional approximations to $\nabla \phi$ as desired. We used this update rule and applied \cref{alg:sweepingScheme} with a $401\times 401$ grid and with approximations to $\nabla \phi$ along the $x$-axis and $y$-axis as well as the $\beta = \pi/4$ direction. The results are seen in \cref{fig:VisibilityExample}, where the yellow set represents the visible set, the black shapes are obstacles and the green dots are the vantage points. In this case, because there is no control variable, the upwind direction is fixed and characteristics are straight lines flowing away from the vantage points. Because of this simple geometry, the scheme requires only one iteration and values at grid nodes are resolved during one of the directional sweeps depending on where they lie relative to the vantage point. For example, if the vantage point is at grid node $(i^*,j^*)$, then the forward-forward sweep will resolve all values $\phi_{ij}$ with $i > i^*$ and $j > j^*$. It should be noted that Oberman and Salvador also devised an upwind sweeping scheme that approximates \eqref{eq:visibilityPDE} with one sweep in each direction by using interpolation to explicitly capture the exact upwind direction. Our method is not an improvement of theirs; we include this example only to demonstrate the diverse applicability of our method. For a full discussion of the visibility problem including rigorous analysis of \eqref{eq:visibilityPDE}, see \cite{ObermanVisibility}.

Our final application is in time-optimal path planning for simple self-driving cars. This problem was first analyzed by Dubins \cite{Dubins} and Reeds and Shepp \cite{ReedsShepp} in a purely geometric sense, and later analyzed in the Hamilton-Jacobi formulation by Takei, Tsai and others \cite{ParkinsonRobot,TakeiTsai2,TakeiTsai1}. Let $(x,y)$ denote the location of the center of mass of the vehicle and $\theta$ denote the orientation. If $W$ is the maximum angular velocity of the car (which enforces a minimum turning radius) and $d$ is the distance from the rear wheels---which drive the car---to the center of mass, then the kinematics are \begin{equation}\label{eq:robotMotion} \begin{split}
\dot x &= v\cos(\theta) - \omega Wd \sin(\theta), \\
\dot y &= v\sin(\theta) + \omega Wd\cos(\theta), \\
\dot \theta &= W\omega,
\end{split} 
\end{equation} where $v,\omega \in [-1,1]$ are normalized control variables representing tangential and angular velocity respectively \cite{Wu}. 

With these kinematics, the optimal travel time function solves the Hamilton-Jacobi equation \begin{equation}\label{eq:HJBRobot}
-1 = \inf_{v,\o} \Big\{[v\cos(\theta) - \o W d \sin(\theta)]\phi_x + [v\sin(\theta) + \o Wd \cos(\theta) ]\phi_y + \o W \phi_\theta \Big\}.
\end{equation} For a full derivation of this equation, we direct the reader to \cite{TakeiTsai2}; they consider the case that $d=0$ so the car is simplified to a point mass, but otherwise the derivation is the same. One notes that the minimization is linear in $(v,\o)$, and thus, since the minimization set $[-1,1] \times [-1,1]$ has finitely many extreme points, there are finitely many values that the pair $(v,\o)$ will take. For technical reasons, one should allow $v \in \{-1,1\}$ and $\o \in \{-1, 0, 1\}$ \cite{TakeiTsai2}.

Equation \eqref{eq:HJBRobot} fits directly into our framework. Discretizing $(x_i,y_j,\theta_k)$, equation \eqref{eq:HJBRobot} is approximated by the update rule \begin{equation}\label{eq:updateRuleRobot} \begin{split}
\phi^*_{ijk}(v,\o) = \Big\{1 + &\frac{\abs{A_k(v,\o)}}{\Delta x} \phi_{i+a_k(u,v),j,k}\\ + &\frac{\abs{B_k(v,\o)}}{\Delta y} \phi_{i,j+b_k(v,\o),k}\\  + &\frac{\abs \o W}{\Delta \theta}\phi_{i,j,k+\text{sign}(\o)}\Big\} \\ &/ \Big\{{\frac{\abs{A_k(v,\o)}}{\Delta x} + \frac{\abs{B_k(v,\o)}}{\Delta y} + \frac{\abs \o W}{\Delta \theta}}\Big\},
\end{split} \end{equation} where \begin{equation}\label{eq:coeffs} \begin{split}
A_k(v,\o) &= v \cos(\theta_k) - \o Wd \sin (\theta_k), \\
B_k(v,\o) &= v \sin (\theta_k) + \o Wd \cos (\theta_k), \\
a_k(v,\o) &= \text{sign}(v \cos (\theta_k) - \o Wd \sin (\theta_k)), \\
b_k(v,\o) &= \text{sign}(v \sin (\theta_k) + \o Wd \cos (\theta_k)).
\end{split}
\end{equation} One can use this update rule in \cref{alg:sweepingScheme} (accounting for three dimensions by performing 8 sweeps per iteration) with the boundary condition $\phi^0_{i^*,j^*,k^*} = 0$ for the desired ending configuration and $\phi^0_{ijk} = +\infty$ otherwise. Then $\phi_{ijk}$ will represent the approximate time needed to travel from grid node $(i,j,k)$ to grid node $(i^*,j^*,k^*)$ while obeying \eqref{eq:robotMotion}.

In three dimensions, it is less obvious how to incorporate grid rotations in a fully principled manner. We discuss this further in \cref{sec:append}. One approach is to restrict ourselves to rotations of the $xy$-plane while keeping the $\theta$-axis fixed. In doing so, we can again trade $(x,y)$ for $(\overline x,\overline y)$ exactly as in the two-dimensional case. Using this strategy, if the rotation angle is $\beta = \arctan(\hat\jmath/\hat\imath)$, the new update rule is \begin{equation}\label{eq:updateRuleRobotRotated} \begin{split}
\overline \phi^*_{ijk}(v,\o) = \Big\{ \Delta s ~+&\abs{\overline A_k(v,\o)} \phi_{i+\hat \imath \overline a_k(u,v),j+\hat \jmath \overline a_k(u,v),k}   \\ +&\abs{\overline B_k(v,\o)} \phi_{i- \hat \imath \overline b_k(v,\o),j+\hat \jmath \overline b_k(v,\o),k}  \\+&\frac{\Delta s\abs \o W}{\Delta \theta}\phi_{i,j,k+\text{sign}(\o)} \Big\} \\ /& \Big\{\abs{A_k(v,\o)} + \abs{B_k(v,\o)} + \frac{\Delta s \abs \o W}{\Delta \theta} \Big\}, \end{split}
\end{equation} where $\Delta s = \sqrt{(\hat \imath \Delta x)^2 + (\hat \jmath \Delta y)^2}$ as before, and \begin{equation}\label{eq:coeffsRot} \begin{split}
\overline A_k(v,\o) &= v \cos (\theta_k+\beta) - \o Wd \sin(\theta_k+\beta), \\
\overline B_k(v,\o) &= v \sin (\theta_k+\beta) + \o Wd \cos (\theta_k+\beta), \\
\overline a_k(v,\o) &= \text{sign}(v \cos (\theta_k+\beta) - \o Wd \sin (\theta_k+\beta)), \\
\overline b_k(v,\o) &= \text{sign}(v \sin (\theta_k+\beta) + \o Wd \cos (\theta_k+\beta)).
\end{split}
\end{equation}

\begin{figure}[b!]
\centering
\begin{subfigure}{0.45\textwidth}
\centering 
\includegraphics[width=\textwidth,trim=0 0 0 0,clip]{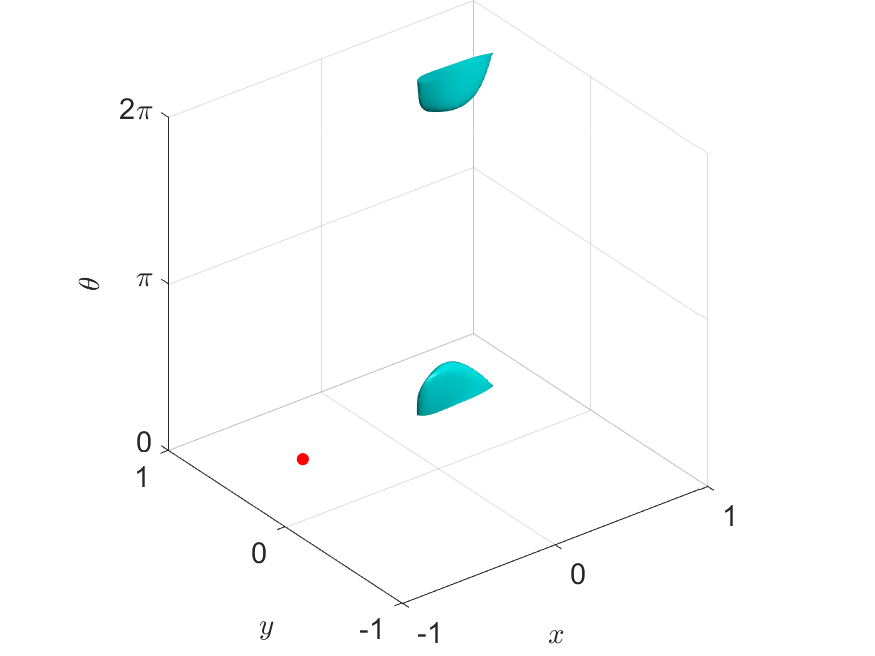}
\caption{Level set $\phi(x,y,\theta) = \tfrac 1 4$.}
\label{fig:valLevel1}
\end{subfigure}~
\begin{subfigure}{0.45\textwidth}
\centering 
\includegraphics[width=\textwidth,trim=0 0 0 0,clip]{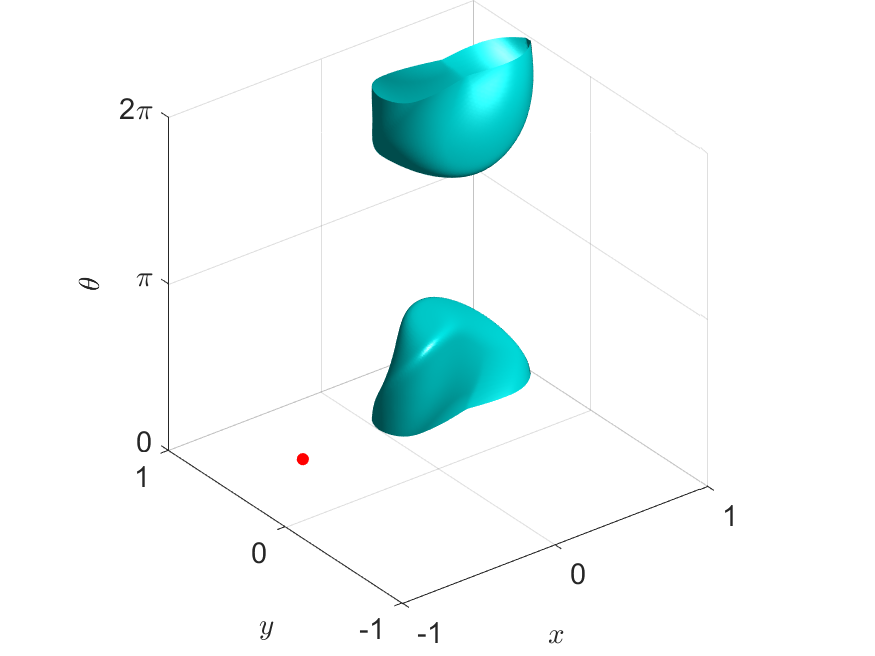}
\caption{Level set $\phi(x,y,\theta) = \tfrac 1 2$.}
\label{fig:valLevel2}
\end{subfigure}\\
\begin{subfigure}{0.45\textwidth}
\centering 
\includegraphics[width=\textwidth,trim=0 0 0 0,clip]{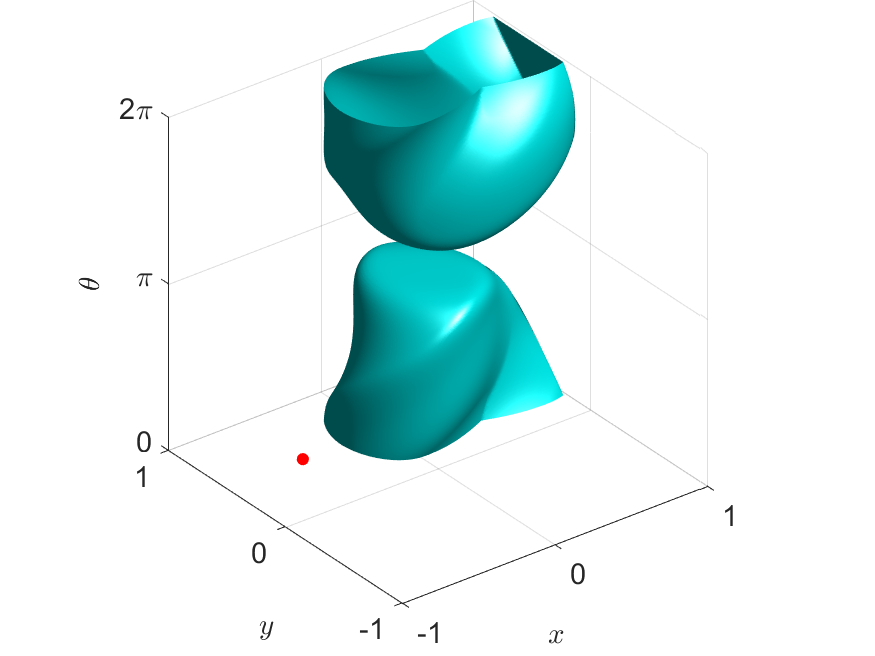}
\caption{Level set $\phi(x,y,\theta) = \tfrac 3 4$.}
\label{fig:valLevel3}
\end{subfigure}~
\begin{subfigure}{0.45\textwidth}
\centering 
\includegraphics[width=\textwidth,trim=0 0 0 0,clip]{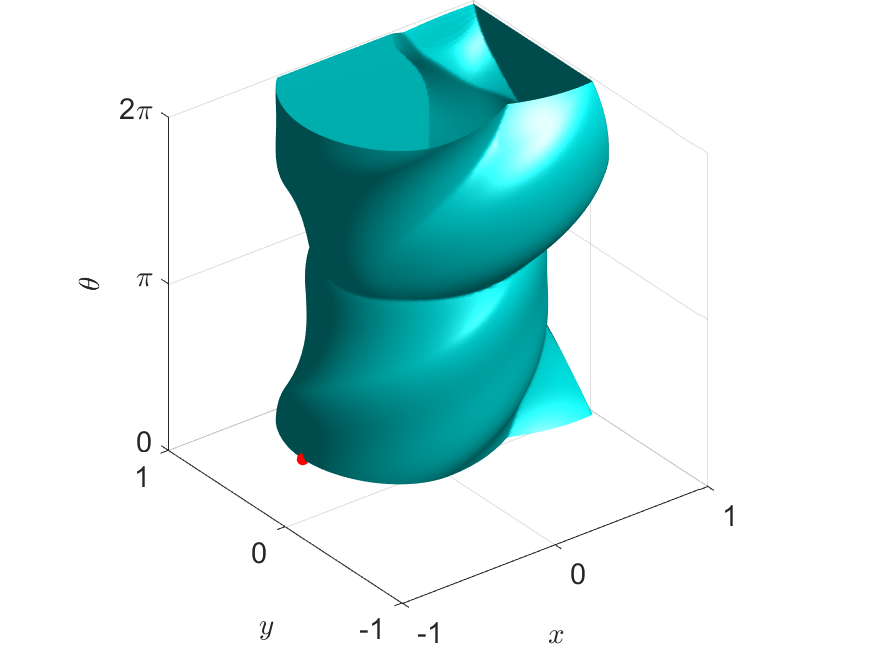}
\caption{Level set $\phi(x,y,\theta) = 1$.}
\label{fig:valLevel4}
\end{subfigure}
\caption{Level sets (cyan) of the travel-time function $\phi(x,y,\theta)$ with ending point $(\tfrac 1 2, \tfrac 1 2, 0)$. Plotted in red is the point $(-\tfrac 1 2, \tfrac 1 2, 0).$ This point should have a travel time of $1$, and indeed the level set $\phi(x,y,\theta)=1$ includes the point.}
\label{fig:valLevelSets}
\end{figure}

\begin{figure}[t!]
\centering
\begin{subfigure}{0.48\textwidth}
\centering
\includegraphics[width=\textwidth,trim=15 40 25 20,clip]{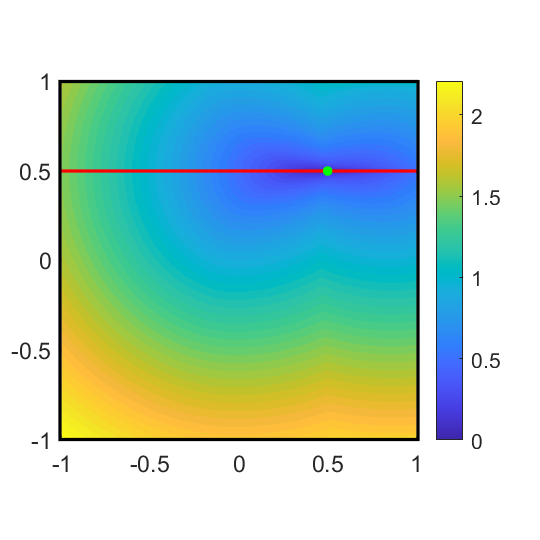}
\caption{Contours of $\phi^{RGU}(x,y,0)$.}
\end{subfigure}~
\begin{subfigure}{0.48\textwidth}
\centering
\includegraphics[width=\textwidth,trim=43 15 35 0,clip]{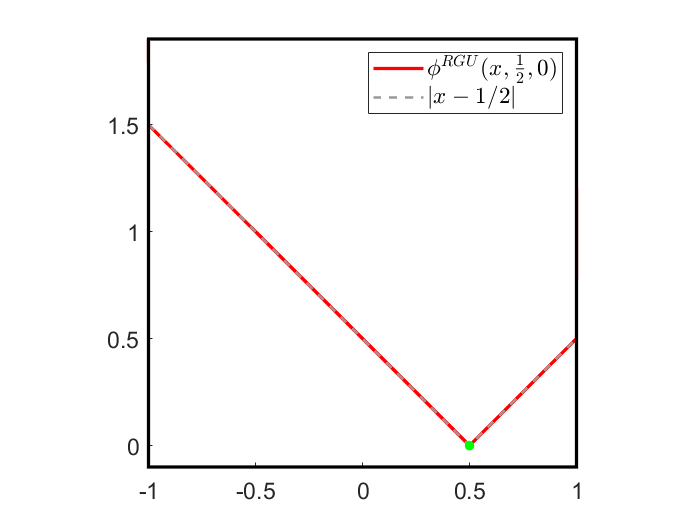}
\caption{Value of $\phi^{RGU}(x,\tfrac 1 2,0)$.}
\end{subfigure}~
\caption{Contour plot of the approximate travel-time function $\phi^{RGU}(x,y,0)$ with ending point $(\tfrac 1 2, \tfrac 1 2, 0)$ [green], computed using our Rotating Grid Upwind (RGU) method. Along the line $(x,\tfrac 1 2,0$) [red] the exact solution value is $\abs{x-1/2}$, and our results match these values.}
\label{fig:robotContourPlot}
\end{figure}

We used these formulas on a $201\times 201\times 201$ discretization of $[-1,1]\times [-1,1]\times [0,2\pi]$ to compute the travel-time function for this control problem when the ending configuration is $(\tfrac 1 2, \tfrac 1 2, 0)$ meaning the car should end at $(x_f,y_f) = (\tfrac 1 2, \tfrac 12)$ facing in the positive $x$-direction. In all these tests, the convergence criterion is $\max_{ijk} \lvert \phi^{n}_{ijk} - \phi^{n-1}\rvert < 10^{-4}$. The results in in \cref{fig:valLevelSets} and \cref{fig:robotContourPlot} were generated using three additional directions to approximate $\phi_x, \phi_y$: the directions of $\beta = \arctan(1/2), \arctan(1), \arctan (2/1)$. One way to evaluate the results is to compare them against known values of the travel-time function. For example, anywhere along the line $(x,\tfrac 1 2, 0)$, the optimal travel time is $\lvert x-\tfrac 1 2\rvert$ since the optimal path simply requires pulling forward or reversing into the final configuration. Accordingly, on the level set plots in \cref{fig:valLevelSets}, we plot the point $(-\tfrac 1 2, \tfrac 1 2,0)$ in red. This point should satisfy $\phi(-\tfrac 1 2, \tfrac 1 2,0) =1$ and indeed, it seems to approximately lie in the level set $\phi(x,y,\theta) = 1$ [\cref{fig:valLevel4}]. Likewise, in \cref{fig:robotContourPlot}, we display the contours of $\phi(x,y,0)$ which show the values of the travel-time function given that the car is facing in the positive $x$-direction. Using these, we can directly compare values of $\phi(x,\tfrac 1 2,0)$ and $\lvert x - \tfrac 1 2\rvert$ and the results line up very well.

\begin{figure}[t!]
\centering
\begin{subfigure}{0.48\textwidth}
\centering
\includegraphics[width=\textwidth,trim=15 40 25 20,clip]{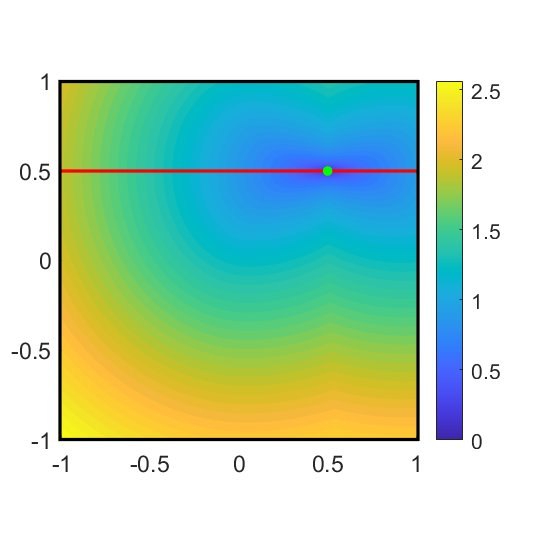}
\caption{Contours of $\phi^{LF}(x,y,0)$.}
\end{subfigure}~
\begin{subfigure}{0.48\textwidth}
\centering
\includegraphics[width=\textwidth,trim=43 15 35 0,clip]{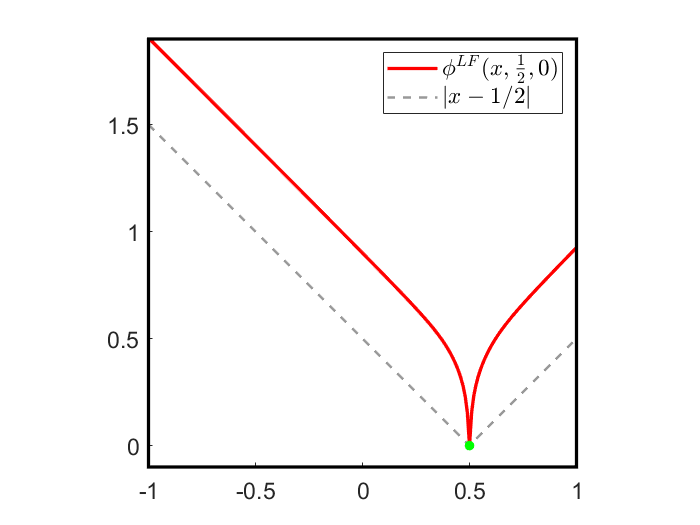}
\caption{Value of $\phi^{LF}(x,\tfrac 1 2,0)$.}
\end{subfigure}~
\caption{Contour plot of the approximate travel-time function $\phi^{LF}(x,y,0)$ with ending point $(\tfrac 1 2, \tfrac 1 2, 0)$ [green], computed using the Lax-Friedrichs (LF) method. Along the line $(x,\tfrac 1 2,0$) [red] the exact solution value is $\abs{x-1/2}$, and due to the diffusive nature of the scheme, the approximate solution incurs some error.}
\label{fig:robotContourPlot_LF}
\end{figure}

Again, we compare our results to those of the Lax-Friedrichs (LF) sweeping scheme \cite{Kao2004}. Because the LF scheme includes artificial viscosity, it has trouble resolving the value function in the neighborhood surrounding the source point $(x_f, y_f, \theta_f)$. Indeed, we computed the solution of the same problem using the LF scheme. Values analogous to those in \cref{fig:robotContourPlot} are displayed in \cref{fig:robotContourPlot_LF}. We note there is some error in the values of $\phi(x,\tfrac 1 2, 0)$. We also notice that the solution suggested by the LF scheme takes larger values throughout the domain, which hints that the optimal travel time is being overestimated. 

\Cref{tab:LFcomp_robot} lists the iteration counts for different grid resolutions, and different solution methods. In the table, ``Basic'' denotes the basic scheme, and ``Basic+$M$''  denotes the basic scheme appended with $M$ grid rotations. A first note is that for this problem, including additional derivative approximations in different directions lowers the number of iterations required for our algorithm to converge. Due to diffusivity, the LF method requires vastly more iterations. In this case, the LF iterations are no more or less efficient than those of our method. Problems where the control values can be resolved explicitly are well-suited to our method. For problems of this type, our method is very likely to ourperform the LF method and is equally easy to implement. It bears repeating that the LF method is more generally applicable and easier to implement for problems with very complicated Hamiltonians \cite{Kao2004}.

\begin{figure}[b!]  
\centering
\begin{subfigure}[t]{0.48\textwidth}
\centering
\includegraphics[width=\textwidth,trim = 50 50 50 30,clip]{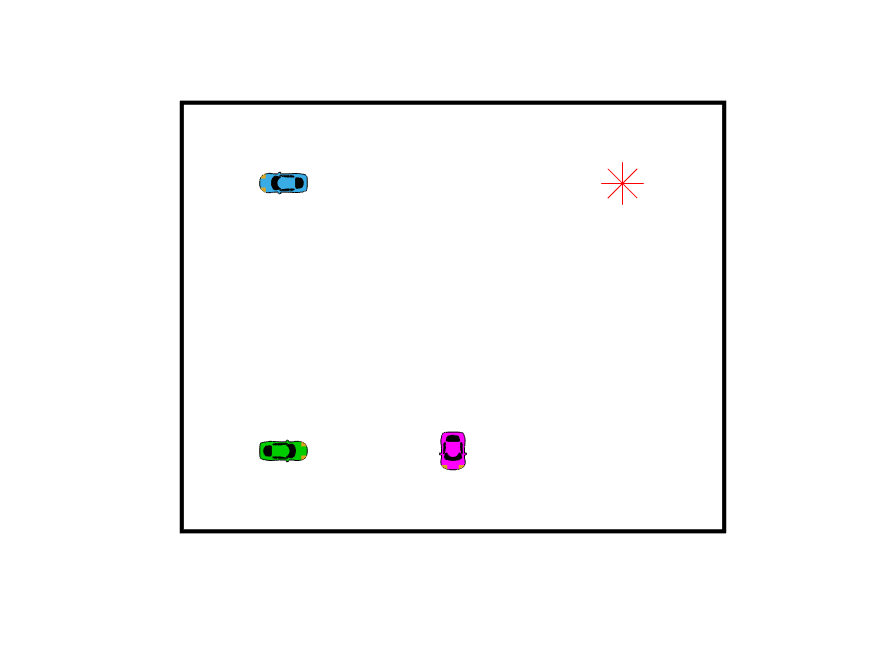}
\caption{Initial configurations.}
\end{subfigure}~ \hfill
\begin{subfigure}[t]{0.48\textwidth}
\centering
\includegraphics[width=\textwidth,trim = 50 50 50 30,clip]{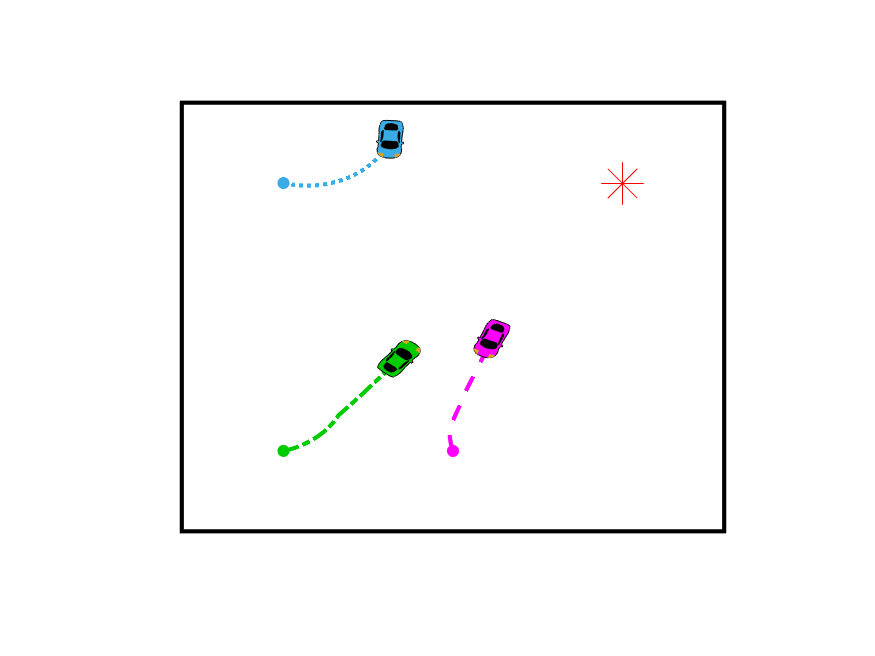}
\caption{$1/3$ of the way along the paths.}
\end{subfigure}\\
\begin{subfigure}[t]{0.48\textwidth}
\centering
\includegraphics[width=\textwidth,trim = 50 50 50 30,clip]{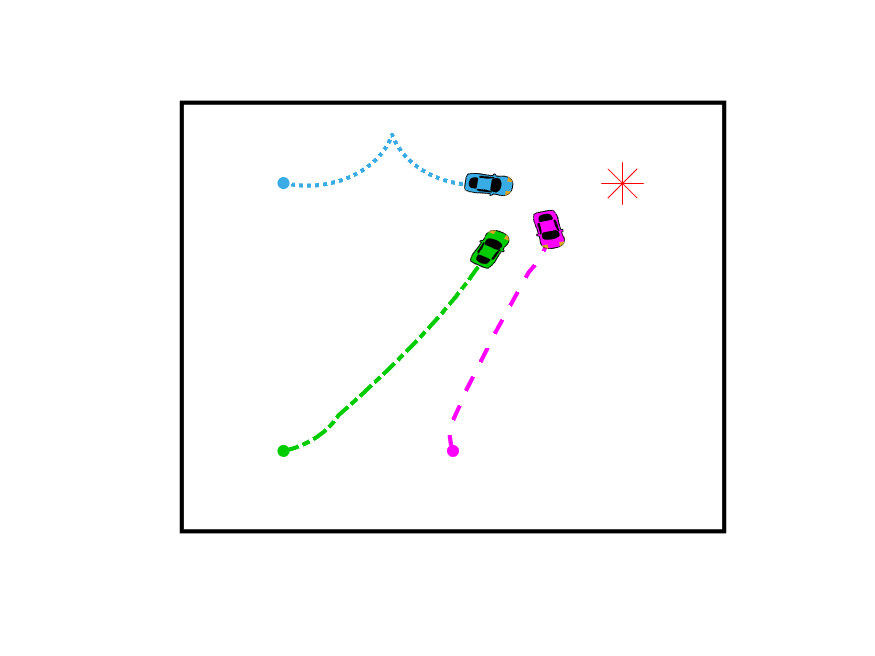}
\caption{$2/3$ of the way along the paths.}
\end{subfigure}~ \hfill
\begin{subfigure}[t]{0.48\textwidth}
\centering
\includegraphics[width=\textwidth,trim = 50 50 50 30,clip]{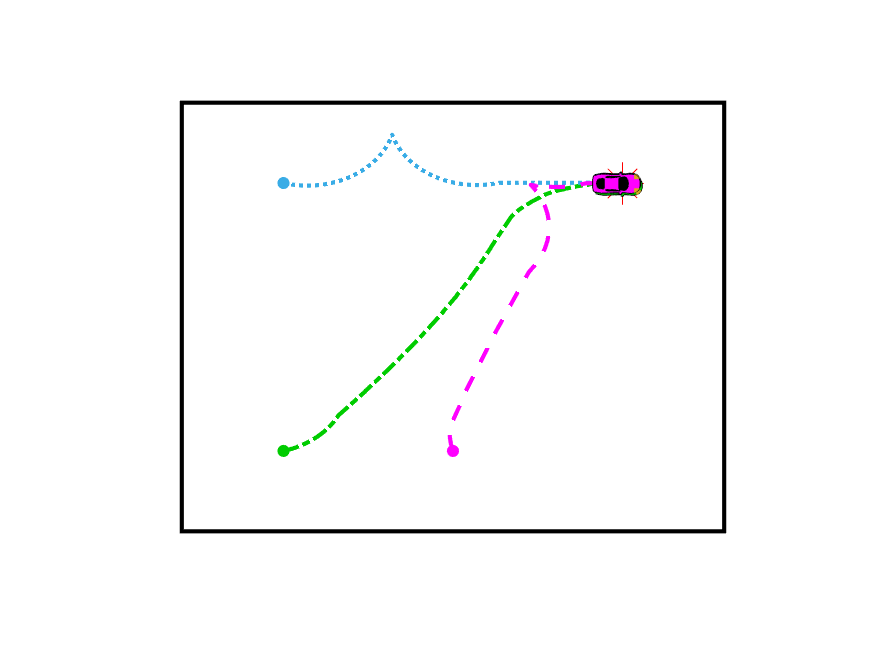}
\caption{Final configurations.}
\end{subfigure}\\
\caption{Optimal paths for cars with initial configurations $(-\tfrac 1 2, \tfrac 1 2, \pi)$ [blue], $(-\tfrac 12, -\tfrac 1 2, 0)$ [green], and $(0, -\tfrac 1 2, \tfrac{5\pi}4)$ [pink]. Final configuration is $(\tfrac 1 2, \tfrac 1 2, 0)$ [red star].}
\label{fig:threePathsFig}
\end{figure}

Another way one can verify the results of these simulations is to compute the actual paths given by the control problem. Having computed the travel-time function $\phi$, one can determine optimal trajectories by integrating \eqref{eq:robotMotion} using control values \begin{equation}\label{eq:controls}
\begin{split} v &= -\text{sign}(\phi_x \cos \theta + \phi_y \sin \theta), \\
\o &= -\text{sign}(-d\phi_x\sin \theta + d\phi_y \cos \theta + \phi_\theta). \end{split} 
\end{equation} Some optimal paths are seen in \cref{fig:threePathsFig}. In those plots, the final location is marked by the red star, and the initial locations are marked by colored dots. The positions of the vehicles are displayed at several points along their respective optimal trajectories. Note, these optimal paths were computed independently and are simply plotted on top of each other; the paths will require different amounts of time to traverse and there is no interaction between the cars. The results appear to agree with a theoretical result of Reeds and Shepp \cite{ReedsShepp} that states that optimal trajectories consist of straight lines and arcs of circles of minimum radius.\\

\begin{table}[t!] 
\centering
{\bf Iteration Counts: Self-Driving Car Example}\\

\begin{tabular}{r | r r r r}
\hline
$I,J, K$  & LF & Basic & Basic$+1$ & Basic$+3$\\
\hline \hline
50 &  99 & 17 & 16 & 17\\
100	& 187 & 25 & 22 & 21\\
200	& 309 & 32 & 26 & 24 \\ \hline
\end{tabular} 
\caption{The iteration counts for different versions of our algorithm, and for the Lax-Friedrichs sweeping scheme, when resolving the optimal travel-time function for the simple self-driving car. Here Basic$+1$ designates the basic method appended with derivative approximations in one additional direction; Basic$+3$ designates the basic method appended with derivative approximations in three additional directions.}
\label{tab:LFcomp_robot}
\end{table}

\section{Conclusion \& Discussion} Fast sweeping methods provide a simple and robust framework for numerical solutions of steady-state Hamilton-Jacobi equations. We have developed a fast sweeping scheme for a class of Hamilton-Jacobi equations arising from steady-state optimal control problems wherein the running cost is independent of the control variables. Our method is exceedingly simple to implement and applies to a wide range of problems. We tested our method against Eikonal equations in different norms, and demonstrated how one can use WENO approximations to improve accuracy. We then suggested a general method for maintaining a square grid, but using approximations to derivatives in rotated directions, so as to more accurately capture the information flow along characteristics. We compare our method against the Lax-Friedrichs method \cite{Kao2004} and demonstrate that in some cases, our method is preferable. Finally, we demonstrated the utility of our method by applying it to two problems arising from engineering applications.    

There are several ways in which our method could be modified or adjusted for other scenarios. We suggest two such modifications now. First, a further exploration of the efficacy of WENO approximations in conjunction with our method could prove interesting. In \cref{sec:WENO}, we demonstrated one method for including WENO approximations, following \cite{ZhangWENOFSM}. However, especially when the solution was non-smooth, we did not achieve the full increase in accuracy that one may desire. It is possible that one could improve this with a closer analysis of the scheme near the point source. One may also try to include WENO approximations with the grid rotations. This is likely to be difficult due to the different sizes of the rotated grid parameters $\Delta \overline x, \Delta \overline y$ which may skew convergence results, so one would need to be cautious. Second, when using a single grid rotation with angle $\beta = \pi/4$, we are essentially using a 9 point stencil for local derivative approximations, which yields a structured triangulation of the domain. It would be interesting to modify the method for unstructured and/or triangulated domains such as those in \cite{Qian1}. In these domains, our method may provide a simpler update rule for Eikonal equations, though a careful analysis would be required. 

\section*{Declarations}

Data sharing not applicable to this article as no datasets were generated or analyzed during the current study. The author has no conflicts of interest to declare that are relevant to the content of this article.

\section*{Acknowledgments}
The author thanks Andrea Bertozzi and Stanley Osher for reading an early version of this manuscript, and for several valuable conversations and suggestions, especially regarding the example of optimal path planning for self-driving cars.

The author also thanks two anonymous reviewers for helpful comments and suggestions which improved the manuscript.

\bibliographystyle{plain}
\bibliography{bibliography}

\begin{thebibliography}{10}

\bibitem{IanMitchell}
K.~Alton and I.~M. Mitchell.
\newblock Optimal path planning under defferent norms in continuous state
  spaces.
\newblock In {\em Proceedings 2006 IEEE International Conference on Robotics
  and Automation, 2006. ICRA 2006.}, pages 866--872, May 2006.

\bibitem{Alton}
Ken Alton and Ian~M Mitchell.
\newblock Fast marching methods for stationary {H}amilton--{J}acobi equations
  with axis-aligned anisotropy.
\newblock {\em SIAM Journal on Numerical Analysis}, 47(1):363--385, 2009.

\bibitem{BardiCapuzzo}
M.~Bardi and I.~Capuzzo-Dolcetta.
\newblock {\em Optimal Control and Viscosity Solutions of
  {H}amilton--{J}acobi--{B}ellman Equations}.
\newblock Modern Birkh{\"a}user Classics. Birkh{\"a}user Boston, 2008.

\bibitem{BarlesSouganidis}
G.~Barles and P.~E. Souganidis.
\newblock Convergence of approximation schemes for fully nonlinear second order
  equations.
\newblock 4:271--283, 1991.
\newblock 3.

\bibitem{Barles}
Guy Barles.
\newblock {\em An Introduction to the Theory of Viscosity Solutions for
  First-Order {H}amilton--{J}acobi Equations and Applications}, pages 49--109.
\newblock Springer Berlin Heidelberg, Berlin, Heidelberg, 2013.

\bibitem{BarlesJakobsen}
Guy. Barles and Espen~R. Jakobsen.
\newblock Error bounds for monotone approximation schemes for
  {H}amilton--{J}acobi--{B}ellman equations.
\newblock {\em SIAM Journal on Numerical Analysis}, 43(2):540--558, 2005.

\bibitem{Bellman}
Richard Bellman.
\newblock The theory of dynamic programming.
\newblock Technical report, {R}and {C}orp, {S}anta {M}onica, {CA}, 1954.

\bibitem{Bellman2}
Richard Bellman.
\newblock {\em Adaptive Control Processes: A Guided Tour}.
\newblock Karreman Mathematics Research Collection, Princeton Legacy Library.
  Princeton University Press, 1961.

\bibitem{BR}
Folkmar Bornemann and Christian Rasch.
\newblock Finite-element discretization of static {H}amilton-{J}acobi equations
  based on a local variational principle.
\newblock {\em Computing and Visualization in Science}, 9(2):57--69, 2006.

\bibitem{Boue}
Michelle Bou\'{e} and Paul Dupuis.
\newblock {M}arkov chain approximations for deterministic control problems with
  affine dynamics and quadratic cost in the control.
\newblock {\em SIAM J. Numer. Anal.}, 36(3):667–695, March 1999.

\bibitem{CaffarelliCrandall}
Luis~A. Caffarelli and Michael~G. Crandall.
\newblock Distance functions and almost global solutions of {E}ikonal
  equations.
\newblock {\em Communications in Partial Differential Equations},
  35(3):391--414, 2010.

\bibitem{YTChow}
Yat~Tin Chow, J{\'e}r{\^o}me Darbon, Stanley Osher, and Wotao Yin.
\newblock Algorithm for overcoming the curse of dimensionality for
  state-dependent {H}amilton-{J}acobi equations.
\newblock {\em Journal of Computational Physics}, 387:376--409, 2019.

\bibitem{CrandallIshiiLions1992}
Michael~G Crandall, Hitoshi Ishii, and Pierre-Louis Lions.
\newblock User's guide to viscosity solutions of second order partial
  differential equations.
\newblock {\em Bulletin of the American mathematical society}, 27(1):1--67,
  1992.

\bibitem{CrandallLions2}
Michael~G Crandall and P-L Lions.
\newblock Two approximations of solutions of {H}amilton-{J}acobi equations.
\newblock {\em Mathematics of computation}, 43(167):1--19, 1984.

\bibitem{CrandallLions}
Michael~G. Crandall and Pierre-Louis Lions.
\newblock Viscosity solutions of {H}amilton-{J}acobi equations.
\newblock {\em Transactions of the American Mathematical Society},
  277(1):1--42, 1983.

\bibitem{RSG3}
Maria~Yuliani Danggo and Sudi Mungkasi.
\newblock A staggered grid finite difference method for solving the elastic
  wave equations.
\newblock {\em Journal of Physics: Conference Series}, 909:012047, Nov 2017.

\bibitem{Darbon}
J{\'e}r{\^o}me Darbon and Stanley Osher.
\newblock Algorithms for overcoming the curse of dimensionality for certain
  {H}amilton--{J}acobi equations arising in control theory and elsewhere.
\newblock {\em Research in the Mathematical Sciences}, 3(1):19, 2016.

\bibitem{Dubins}
L.~E. Dubins.
\newblock On curves of minimal length with a constraint on average curvature,
  and with prescribed initial and terminal positions and tangents.
\newblock {\em American Journal of Mathematics}, 79(3):497--516, 1957.

\bibitem{Engquist1}
Bj{\"o}rn Engquist, Brittany~D Froese, and Yen-Hsi~Richard Tsai.
\newblock Fast sweeping methods for hyperbolic systems of conservation laws at
  steady state.
\newblock {\em Journal of Computational Physics}, 255:316--338, 2013.

\bibitem{Engquist2}
Bj{\"o}rn Engquist, Brittany~D Froese, and Yen-Hsi~Richard Tsai.
\newblock Fast sweeping methods for hyperbolic systems of conservation laws at
  steady state {II}.
\newblock {\em Journal of Computational Physics}, 286:70--86, 2015.

\bibitem{EvansControlNotes}
Lawrence~C Evans.
\newblock An introduction to mathematical optimal control theory version 0.2.
\newblock {\em Lecture notes available online}.

\bibitem{RSG2}
Kai Gao and Lianjie Huang.
\newblock An improved rotated staggered-grid finite-difference method with
  fourth-order temporal accuracy for elastic-wave modeling in anisotropic
  media.
\newblock {\em Journal of Computational Physics}, 350:361 -- 386, 2017.

\bibitem{JiangWENO}
Guang-Shan Jiang and Danping Peng.
\newblock Weighted {ENO} schemes for {H}amilton--{J}acobi equations.
\newblock {\em SIAM Journal on Scientific computing}, 21(6):2126--2143, 2000.

\bibitem{Kao20042}
Chiu~Yen Kao, Carmeliza Navasca, and Stanley Osher.
\newblock The {L}ax-{F}riedrichs sweeping method for optimal control problems
  in continuous and hybrid dynamics.
\newblock {\em Nonlinear Analysis: Theory, Methods \& Applications}, 63(5):1561
  -- 1572, 2005.
\newblock Invited Talks from the Fourth World Congress of Nonlinear Analysts
  (WCNA 2004).

\bibitem{Kao2004}
Chiu~Yen Kao, Stanley Osher, and Jianliang Qian.
\newblock {L}ax--{F}riedrichs sweeping scheme for static {H}amilton--{J}acobi
  equations.
\newblock {\em Journal of Computational Physics}, 196(1):367--391, 2004.

\bibitem{Kao2005}
Chiu-Yen. Kao, Stanley. Osher, and Yen-Hsi. Tsai.
\newblock Fast sweeping methods for static {H}amilton--{J}acobi equations.
\newblock {\em SIAM Journal on Numerical Analysis}, 42(6):2612--2632, 2005.

\bibitem{Kao}
Chiu-Yen Kao and Richard Tsai.
\newblock Properties of a level set algorithm for the visibility problems.
\newblock {\em Journal of Scientific Computing}, 35:170--191, 2008.

\bibitem{LinCurse}
Alex~Tong Lin, Yat~Tin Chow, and Stanley~J. Osher.
\newblock A splitting method for overcoming the curse of dimensionality in
  {H}amilton--{J}acobi equations arising from nonlinear optimal control and
  differential games with applications to trajectory generation.
\newblock {\em Communications in Mathematical Sciences}, 16(7), 1 2018.

\bibitem{Luo3}
Songting Luo.
\newblock A uniformly second order fast sweeping method for {E}ikonal
  equations.
\newblock {\em Journal of Computational Physics}, 241:104--117, 2013.

\bibitem{Luo1}
Songting Luo, Shingyu Leung, and Jianliang Qian.
\newblock An adjoint state method for numerical approximation of continuous
  traffic congestion equilibria.
\newblock {\em Communications in Computational Physics}, 10, 11 2011.

\bibitem{Luo4}
Songting Luo, Jianliang Qian, and Robert Burridge.
\newblock High-order factorization based high-order hybrid fast sweeping
  methods for point-source {E}ikonal equations.
\newblock {\em SIAM Journal on Numerical Analysis}, 52(1):23--44, 2014.

\bibitem{Luo2}
Songting Luo, Jianliang Qian, and Plamen Stefanov.
\newblock Adjoint state method for the identification problem in {SPECT}:
  Recovery of both the source and the attenuation in the attenuated x-ray
  transform.
\newblock {\em SIAM Journal on Imaging Sciences}, 7(2):696--715, 2014.

\bibitem{LuoZhao}
Songting Luo and Hongkai Zhao.
\newblock Convergence analysis of the fast sweeping method for static convex
  {H}amilton--{J}acobi equations.
\newblock {\em Research in the Mathematical Sciences}, 3(1):35, 2016.

\bibitem{ObermanVisibility}
Adam Oberman and Tiago Salvador.
\newblock A partial differential equation obstacle problem for the level set
  approach to visibility.
\newblock {\em Journal of Scientific Computing}, 82(1):14, 2020.

\bibitem{ObermanVladimirsky}
Adam~M. Oberman, Ryo Takei, and Alexander Vladimirsky.
\newblock Homogenization of metric {H}amilton--{J}acobi equations.
\newblock {\em Multiscale Modeling \& Simulation}, 8(1):269--295, 2009.

\bibitem{Osher1993}
Stanley Osher.
\newblock A level set formulation for the solution of the {D}irichlet problem
  for {H}amilton-{J}acobi equations.
\newblock {\em SIAM Journal on Mathematical Analysis}, 24(5):1145--1152, 1993.

\bibitem{OsherFedkiw}
Stanley Osher and Ronald~P. Fedkiw.
\newblock {\em Level set methods and dynamic implicit surfaces}, volume 153 of
  {\em Applied Mathematical Sciences}.
\newblock Springer--Verlag, 2003.

\bibitem{OsherShu}
Stanley Osher and Chi-Wang Shu.
\newblock High order essentially non--oscillatory schemes for
  {H}amilton--{J}acobi equations.
\newblock {\em SIAM Journal of Numerical Analysis}, 28(4):907--922, August
  1991.

\bibitem{OsterWilson}
George~F Oster and Edward~O Wilson.
\newblock {\em Caste and ecology in the social insects}.
\newblock Princeton University Press, 1978.

\bibitem{ParkinsonRobot}
C.~{Parkinson}, A.~L. {Bertozzi}, and S.~J. {Osher}.
\newblock A {H}amilton-{J}acobi formulation for time-optimal paths of
  rectangular nonholonomic vehicles.
\newblock In {\em 2020 59th IEEE Conference on Decision and Control (CDC)},
  pages 4073--4078, 2020.

\bibitem{Parkinson}
Christian Parkinson, David Arnold, Andrea~L Bertozzi, Yat~Tin Chow, and Stanley
  Osher.
\newblock Optimal human navigation in steep terrain: a
  {H}amilton--{J}acobi--{B}ellman approach.
\newblock {\em Communications in Mathematical Sciences}, 17(1):227--242, 2019.

\bibitem{Pham}
Huy\^en Pham.
\newblock {\em Continuous-time Stochastic Optimal Control and Optimization with
  Financial Applications}.
\newblock Springer-Verlag Berlin Heidelberg, 1 edition, 2009.

\bibitem{Qian2}
Jianliang Qian, Yong-Tao Zhang, and Hong-Kai Zhao.
\newblock A fast sweeping method for static convex {H}amilton--{J}acobi
  equations.
\newblock {\em Journal of Scientific Computing}, 31(1-2):237--271, 2007.

\bibitem{Qian1}
Jianliang Qian, Yong-Tao Zhang, and Hong-Kai Zhao.
\newblock Fast sweeping methods for {E}ikonal equations on triangular meshes.
\newblock {\em SIAM Journal on Numerical Analysis}, 45(1):83--107, 2007.

\bibitem{RSG4}
Li~Qin, Ma~{Sui-Bo}, Zhao Bin, and Zhang Wei.
\newblock An improved rotated staggered grid finite difference scheme in coal
  seam.
\newblock {\em Applied Geophysics}, 2019.

\bibitem{ReedsShepp}
J.~A. Reeds and L.~A. Shepp.
\newblock Optimal paths for a car that goes both forwards and backwards.
\newblock {\em Pacific J. Math.}, 145(2):367--393, 1990.

\bibitem{RSG}
Erik~H. Saenger, Norbert Gold, and Serge~A. Shapiro.
\newblock Modeling the propagation of elastic waves using a modified
  finite-difference grid.
\newblock {\em Wave Motion}, 31(1):77 -- 92, 2000.

\bibitem{Sethian1}
J~A Sethian.
\newblock A fast marching level set method for monotonically advancing fronts.
\newblock {\em Proceedings of the National Academy of Sciences},
  93(4):1591--1595, 1996.

\bibitem{SethVlad2}
James~A. Sethian and A.~Vladimirsky.
\newblock Ordered upwind methods for static {H}amilton-{J}acobi equations.
\newblock {\em Proceedings of the National Academy of Sciences},
  98(20):11069--11074, 2001.

\bibitem{SethVlad1}
James~A. Sethian and A.~Vladimirsky.
\newblock Ordered upwind methods for static {H}amilton-{J}acobi equations:
  Theory and algorithms.
\newblock {\em SIAM Journal on Numerical Analysis}, 41(1):325--363, 2003.

\bibitem{Shu}
Chi-Wang Shu.
\newblock High order numerical methods for time dependent {H}amilton--{J}acobi
  equations.
\newblock In {\em Mathematics and computation in imaging science and
  information processing}, pages 47--91. World Scientific, 2007.

\bibitem{BangBang}
LM~Sonneborn and FS~{Van Vleck}.
\newblock The bang-bang principle for linear control systems.
\newblock {\em Journal of the Society for Industrial and Applied Mathematics,
  Series A: Control}, 2(2):151--159, 1964.

\bibitem{Souganidis1985}
Panagiotis~E Souganidis.
\newblock Approximation schemes for viscosity solutions of {H}amilton-{J}acobi
  equations.
\newblock {\em Journal of Differential Equations}, 59(1):1 -- 43, 1985.

\bibitem{TakeiTsai1}
R.~Takei, R.~Tsai, H.~Shen, and Y.~Landa.
\newblock A practical path-planning algorithm for a simple car: a
  {H}amilton-{J}acobi approach.
\newblock In {\em Proceedings of the 2010 American Control Conference}, pages
  6175--6180, June 2010.

\bibitem{TakeiTsai2}
Ryo Takei and Richard Tsai.
\newblock Optimal trajectories of curvature constrained motion in the
  {H}amilton-{J}acobi formulation.
\newblock {\em Journal of Scientific Computing}, 54(2):622--644, Feb 2013.

\bibitem{Tsai}
Y.-H.R. Tsai, L.-T. Cheng, S.~Osher, P.~Burchard, and G.~Sapiro.
\newblock Visibility and its dynamics in a {PDE} based implicit framework.
\newblock {\em Journal of Computational Physics}, 199(1):260 -- 290, 2004.

\bibitem{TsaiOsherSweep}
Yen-Hsi~Richard. Tsai, Li-Tien. Cheng, Stanley. Osher, and Hong-Kai. Zhao.
\newblock Fast sweeping algorithms for a class of {H}amilton--{J}acobi
  equations.
\newblock {\em SIAM Journal on Numerical Analysis}, 41(2):673--694, 2003.

\bibitem{Tsitsiklis}
J.~N. Tsitsiklis.
\newblock Efficient algorithms for globally optimal trajectories.
\newblock {\em IEEE Transactions on Automatic Control}, 40(9):1528--1538, Sep
  1995.

\bibitem{RSG5}
Kang Wang, Suping Peng, Yongxu Lu, and Xiaoqin Cui.
\newblock The velocity-stress finite-difference method with a rotated staggered
  grid applied to seismic wave propagation in a fractured medium.
\newblock {\em Geophysics}, 85(2):T89--T100, 2020.

\bibitem{Wu}
Weiguo Wu, Huitang Chen, and Peng-Yung Woo.
\newblock Time optimal path planning for a wheeled mobile robot.
\newblock {\em Journal of Robotic Systems}, 17(11):585--591, 2000.

\bibitem{RSG1}
Lei {Yang}, Hongyong {Yan}, and Hong {Liu}.
\newblock {Optimal rotated staggered-grid finite-difference schemes for elastic
  wave modeling in {TTI} media}.
\newblock {\em Journal of Applied Geophysics}, 122:40--52, November 2015.

\bibitem{ZhangWENOFSM}
Yong-Tao Zhang, Hong-Kai Zhao, and Jianliang Qian.
\newblock High order fast sweeping methods for static {H}amilton--{J}acobi
  equations.
\newblock {\em Journal of Scientific Computing}, 29(1):25--56, 2006.

\bibitem{ZhaoOsher}
Hong-Kai Zhao, Stanley Osher, Barry Merriman, and Myungjoo Kang.
\newblock Implicit and nonparametric shape reconstruction from unorganized data
  using a variational level set method.
\newblock {\em Computer Vision and Image Understanding}, 80(3):295 -- 314,
  2000.

\bibitem{Zhao}
Hongkai Zhao.
\newblock A fast sweeping method for {E}ikonal equations.
\newblock {\em Mathematics of computation}, 74(250):603--627, 2005.

\bibitem{TomlinReachAvoid}
Z.~Zhou, J.~Ding, H.~Huang, R.~Takei, and C.~Tomlin.
\newblock Efficient path planning algorithms in reach-avoid problems.
\newblock {\em Automatica}, 89:28 -- 36, 2018.

\end{thebibliography}

\appendix
\section{Appendix: 3D Implementation} \label{sec:append}

In this appendix, we briefly describe the implementation of the rotating-grid method in three dimensions. In this case, the equation of interest is \begin{equation} \label{eq:3deqn}
-r(x,y,z) = \inf_{a \in A} \bigg\{ f_1(x,y,z,a)\phi_x + f_2(x,y,z,a) \phi_y + f_3(x,y,z,a) \phi_z \bigg\}.
\end{equation} The extension of the basic method to 3D is straightforward. We discretize the domain into $(x_i,y_j,z_k)$ with uniform grid parameters $\Delta x, \Delta y$ and $\Delta z$. Then following the work in \cref{sec:ourBasicScheme}, we arrive at the local upwind approximation \begin{equation} \label{eq:3dupdate}
\phi^*_{ijk}(a) = \frac{r_{ijk} + 
\frac{\abs{f_{1,ijk}(a)}}{\Delta x} \phi_{i+\xi_{1,ijk}(a),j,k} + \frac{\abs{f_{2,ijk}(a)}}{\Delta y} \phi_{i,j+\xi_{2,ijk}(a),k} + \frac{\abs{f_{3,ijk}(a)}}{\Delta z} \phi_{i,j,k+\xi_{3,ijk}(a)}  
}{\frac{\abs{f_{1,ijk}(a)}}{\Delta x} + \frac{\abs{f_{2,ijk}(a)}}{\Delta y} + \frac{\abs{f_{3,ijk}(a)}}{\Delta z}},
\end{equation} where $\xi_{\ell,ijk}(a) = \text{sign}(f_{\ell,ijk}(a))$ as before. One can then use the update rule $\phi^{n}_{ijk} = \min\{ \phi^{n-1}_{ijk}, \min_{a\in A} \phi^*_{ijk}(a) \}$, while performing 8 sweeps per iteration to account for all combinations of sweeping directions. 

In theory, introducing a rotation in $\R^3$ is not too different from introducing a rotation in $\R^2$. One can choose an orthogonal matrix $U = [\,\, u_1 \,\, | \,\, u_2 \,\, | \,\, u_3\,\,] $ whose columns represent the directions of the new axes, and set \begin{equation} \begin{pmatrix} \overline x \\ \overline y \\ \overline z \end{pmatrix} = U^t \begin{pmatrix} x \\ y \\  z \end{pmatrix} \,\,\,\,\,\,\,\,\, \text { so that } \,\,\,\,\,\,\,\,\,\,\nabla \phi = U \overline \nabla \phi, \end{equation} where $\nabla$ represents the gradient in the original coordinates, and $\overline \nabla $ represents the gradient with respect to the new coordinates. Plugging this representation of $\nabla \phi$ into \eqref{eq:3deqn} and denoting $ f \defeq (f_1, f_2, f_3)$ gives \begin{equation}  \label{eq:3drotatedeqngrad}
-r = \inf_{a \in A}  \bigg \{ \left\langle f(a), U \overline \nabla \phi \right \rangle \bigg\} = \inf_{a \in A} \bigg\{ \left\langle U^t  f(a) , \overline \nabla \phi \right\rangle \bigg\}.
\end{equation} Thus, defining $\overline f_\ell(\overline x, \overline y, \overline z,a) =\langle u_\ell , f(\overline x, \overline y, \overline z,a) \rangle$ for $\ell = 1,2,3$, the rotated equation is \begin{equation} \label{eq:3drotatedeqn}
-r(\overline x, \overline y, \overline z) = \inf_{a \in A} \bigg\{ \overline f_1 (\overline x, \overline y, \overline z,a) \phi_{\overline x} + \overline f_2 (\overline x, \overline y, \overline z,a) \phi_{\overline y} + \overline f_3 (\overline x, \overline y, \overline z,a) \phi_{\overline z}\bigg\}.
\end{equation} 

Now the question arises of how to discretize this equation. Similar to the 2D formulation, we would like to avoid defining a new grid, but rather restrict ourselves to rotations which allow us to use the already-defined grid points to approximate derivatives in different directions. There is a practical complication to address here. As demonstrated above, in 2D it is sufficient to choose a grid point $(\hat \imath, \hat \jmath)$ and rotate the grid so that the $\overline x$-axis points at $(\hat \imath, \hat \jmath)$. Having done so, the new $\overline y$-axis points at $(-\hat \jmath, \hat \imath)$ as shown in \cref{fig:gridRotate}. However, in 3D, there are infinitely many rotations which fix the $\overline x$-axis in a specified direction. Thus, in analogy to the 2D scenario, there needs to be a principled manner by which to point the $\overline x$-axis toward a desired grid point $(\hat \imath, \hat \jmath, \hat k)$ while ensuring that the $\overline y$- and $\overline z$-axes are still pointed toward other grid points, so as to avert the need to define a new grid. 

We suggest two ways for doing this. The first, which is simpler but not as general, is to restrict oneself to rotations which fix one of the axes, as we did in the example of time-optimal path planning for self-driving cars in \cref{sec:otherApps}. Here one chooses $(\hat \imath, \hat \jmath)$ as before, and also specifies which axis is to remain fixed.  In doing so, the 3D implementation is effectively reduced to a 2D implementation, since one of the derivatives follows through the computation without changing. 

\begin{figure}[h!]
\centering
\includegraphics[width=0.30\textwidth]{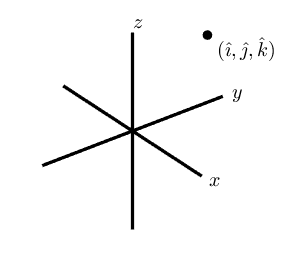}\,\,\, \includegraphics[width=0.30\textwidth]{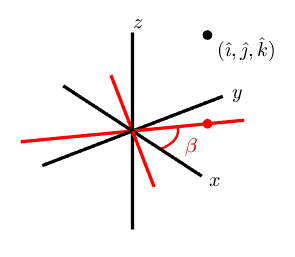} \,\,\, \includegraphics[width=0.30\textwidth]{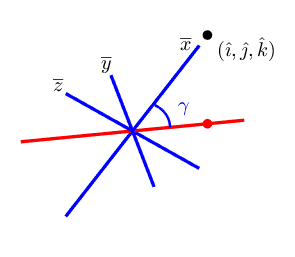} \\
\caption{We orient $\overline x$ to point at $(\hat \imath, \hat \jmath, \hat k)$ by viewing it as the image of the $x$-axis under successive rotations: first a rotation by $\beta$ about the $z$-axis, and then a rotation by $\gamma$ about the line $\hat \imath x + \hat \jmath y = 0$. }
\label{fig:3drotfig}
\end{figure}

The second method can handle general rotations, but is slightly more difficult to describe. Here, we suggest choosing a grid point $(\hat \imath, \hat \jmath, \hat k)$ and using the rotation which orients the $\overline x$-axis toward $(\hat \imath, \hat \jmath, \hat k)$ by viewing it as the image of the $x$-axis under two successive rotations: first, a rotation by $\beta = \arctan(\hat \jmath / \hat \imath)$ about the $z$-axis, and then a rotation by $\gamma = \arctan(\hat k /  \sqrt{\hat \imath ^2 + \hat \jmath ^2})$ about the line $\hat \imath x + \hat \jmath y = 0$. This is illustrated in \cref{fig:3drotfig} where the black lines are the original axes, the red lines are the axes resulting from the first rotation (note, the $z$-axis is unchanged under the first rotation), and the blue lines are the new axes after both rotations. 

In terms of the old coordinates, the new orthogonal coordinates are given by \begin{equation} \label{eq:3dDirections} \begin{split} &\text{positive } \overline x \text{ direction:} \,\ (\hat \imath, \hat \jmath, \hat k),\\ &\text{positive } \overline y \text{ direction:} \,\, (-\hat \jmath, \hat \imath, 0), \\ &\text{positive } \overline z \text{ direction:} \,\, (-\hat \imath \hat k, -\hat \jmath \hat k, \hat \imath^2 + \hat \jmath^2). \end{split} \end{equation} As formulated in \eqref{eq:3drotatedeqngrad}, the columns of $U$ are normalized versions of the three vectors in \eqref{eq:3dDirections}. Using these, we can write the upwind approximations to the derivatives $\phi_{\overline x}, \phi_{\overline y}$ and $\phi_{\overline z}$ necessary to approximate \eqref{eq:3drotatedeqn}. Indeed, \begin{equation} \label{eq:upwindApproxRot3d} \begin{split}
\left (\overline f_1(\overline x, \overline y, \overline z,a)\phi_{\overline x}\right)_{ijk} &= \abs{\overline f_{1,ijk}(a)} \frac{\phi_{i+\hat \imath \overline \xi_{1,ijk}(a),j+ \hat \jmath \overline \xi_{1,ijk}(a), k + \hat k \overline \xi_{1,ijk}(a)} - \phi_{ijk}}{\Delta \overline x}, \\ & \\ \left (\overline f_2(\overline x, \overline y, \overline z,a)\phi_{\overline y}\right)_{ijk} &= \abs{\overline f_{2,ijk}(a)} \frac{\phi_{i-\hat \jmath \overline \xi_{2,ijk}(a),j+ \hat \imath \overline \xi_{2,ijk}(a), k} - \phi_{ijk}}{\Delta \overline y}, \\ & \\ \left (\overline f_3(\overline x, \overline y, \overline z,a)\phi_{\overline z}\right)_{ijk} &= \abs{\overline f_{3,ijk}(a)} \frac{\phi_{i-\hat \imath \hat k \overline \xi_{3,ijk}(a),j- \hat \jmath\hat k \overline \xi_{3,ijk}(a), k + (\hat \imath^2 + \hat \jmath^2)\overline \xi_{3,ijk}(a)} - \phi_{ijk}}{\Delta \overline z}
\end{split} \end{equation} where, as before, $\overline \xi_{\ell,ijk}(a) = \text{sign}\left(\overline f_{\ell,ijk}(a)\right),$ and the new grid parameters are given by \begin{equation} \label{eq:3dgridparams}
\begin{split}
\Delta \overline x &= \sqrt{(\hat \imath \Delta x)^2 + (\hat \jmath \Delta y)^2 + (\hat k \Delta z)^2}, \\ 
\Delta \overline y &= \sqrt{(\hat \jmath \Delta x)^2 + (\hat \imath \Delta y)^2 }, \\
\Delta \overline x &= \sqrt{(\hat \imath \hat k \Delta x)^2 + (\hat \jmath \hat k \Delta y)^2 + ((\hat \imath^2 + \hat \jmath^2)\Delta z)^2}. 
\end{split}
\end{equation}

We will demonstrate both of these methods of implementation using the 1-norm Eikonal equation as an example; first for its simplicity, and second because the level sets of the solution have sharp edges which allow us to easily verify the results. In 3D, the 1-norm Eikonal equation is given by \begin{equation} \label{eq:1normeik3d} -1 = -\abs{\phi_x} - \abs{\phi_y} - \abs{\phi_z} = \inf_{a_i \in \{\pm 1\}} \bigg\{a_1 \phi_x + a_2 \phi_y + a_3 \phi_z \bigg\}.\end{equation} We use the boundary condition $\phi(0,0,0) = 0$. The solution is $\phi(x,y,z) = \max\{\abs x, \abs y, \abs z\}$, and the level sets of this solution are perfect cubes.  In all cases, we use a $201 \times 201 \times 201$ discretization of $[-1,1]^3$ and display the level set $\phi(x,y,z) = 1/2$ which should be a cube centered at the origin of side length 1.  We successively build better approximations of the solution by including additional approximations to the derivatives in different rotated directions.

\begin{figure}[h!]
\centering
\begin{subfigure}{0.45\textwidth}
\centering 
\includegraphics[width=\textwidth]{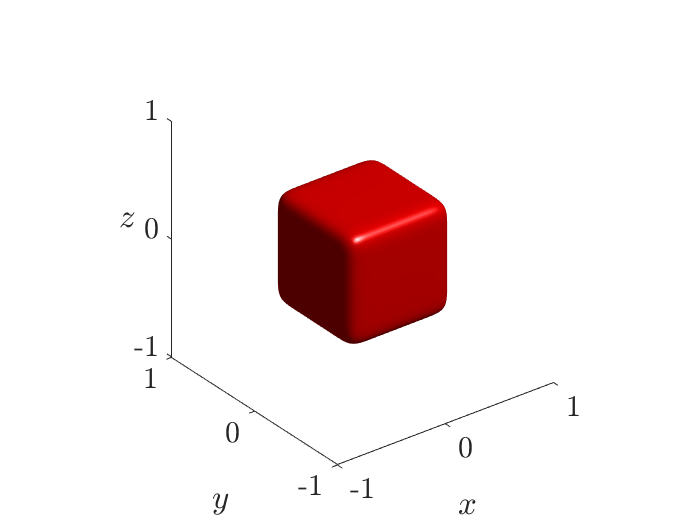} 
\caption{Level set of the solution using the update rule \eqref{eq:updaterule3deik1basic}.}
\label{fig:lastsubfig1}
\end{subfigure} \,\, 
\begin{subfigure}{0.45\textwidth}
\centering 
\includegraphics[width=\textwidth]{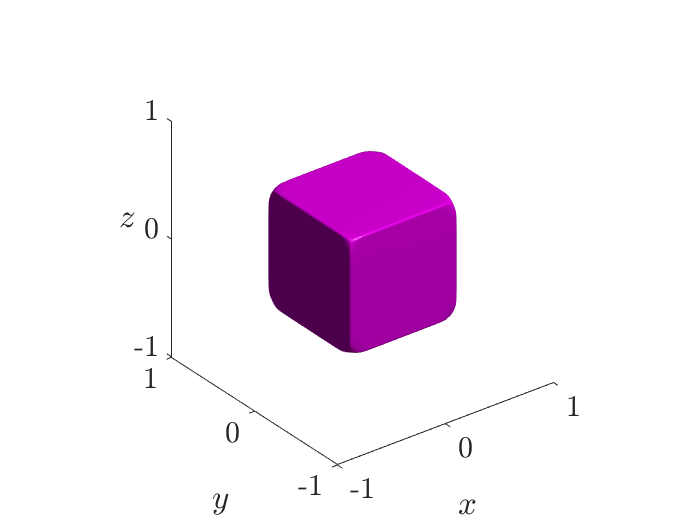} 
\caption{Level set of the solution using the update rule \eqref{eq:updaterule3deik1edges}.}
\label{fig:lastsubfig2}
\end{subfigure} \\
\begin{subfigure}{0.45\textwidth}
\centering 
\includegraphics[width=\textwidth]{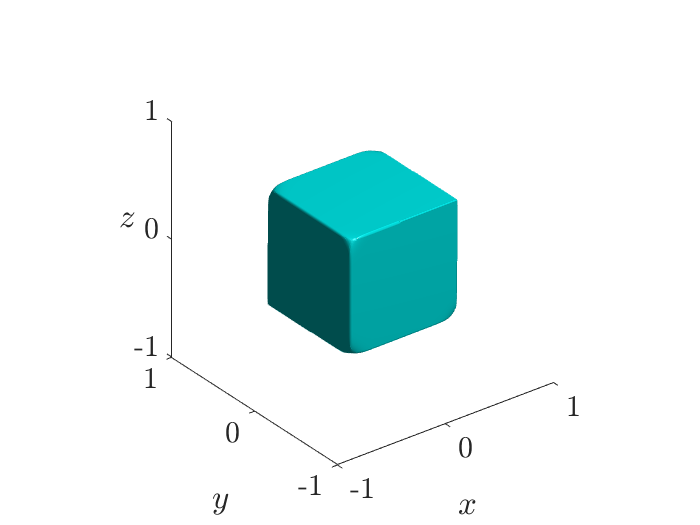} 
\caption{Level set of the solution using the update rule \eqref{eq:updateRule1Corner}.}
\label{fig:lastsubfig3}
\end{subfigure}\,\,
\begin{subfigure}{0.45\textwidth}
\centering 
\includegraphics[width=\textwidth]{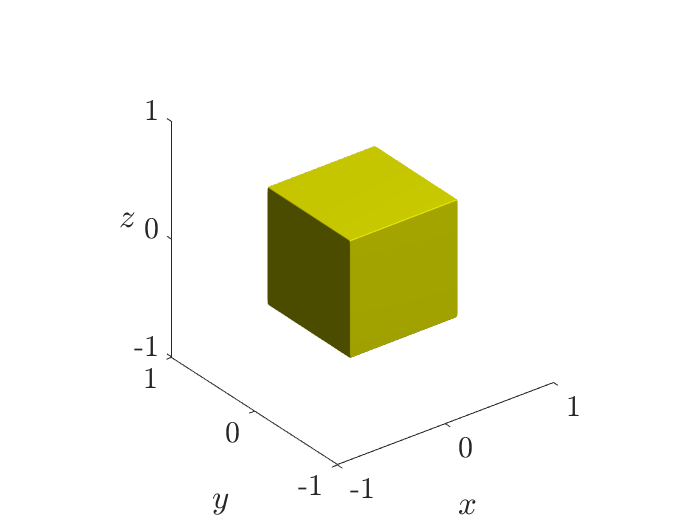} 
\caption{Level set of the solution using the update rule \eqref{eq:updateRuleAllCorners}.}
\label{fig:lastsubfig4}
\end{subfigure}
\caption{The $1/2$-level set of the approximate solution to $\|\nabla \phi\|_1 =1$ resulting from our method when using the update rules  \eqref{eq:updaterule3deik1basic} (no rotated directions), \eqref{eq:updaterule3deik1edges} (rotations which keep one axis fixed), \eqref{eq:updateRule1Corner} (an additional rotation to resolve the corners along $(1,-1,1)$), and \eqref{eq:updateRuleAllCorners} (rotations to capture all the corners).}
\label{fig:eik1norm3d}
\end{figure}

 Applying the basic method, we find that the local upwind approximation to the solution is \begin{equation} \label{eq:updaterule3deik} \phi_{ijk}^*(a) = \frac{1 + \frac 1 {\Delta  x} \phi_{i+\text{sign}(a_1),j,k} + \frac 1 {\Delta  y} \phi_{i,j+\text{sign}(a_2),k} + \frac 1 {\Delta z} \phi_{i,j,k+\text{sign}(a_3)}}{\frac 1 {\Delta x} + \frac 1 {\Delta y} + \frac 1 {\Delta z}}.\end{equation} Then the local update rule for the iteration is \begin{equation} \label{eq:updaterule3deik1basic} \phi^n_{ijk} = \min\bigg\{\phi^{n-1}_{ijk}, \,\,\,\min_{a}  \phi^*_{ijk}(a) \bigg\}. \end{equation} A level set of the solution produced by this update rule is displayed in red in \cref{fig:lastsubfig1}. Note the rounding along the edges and at the corners. 

Next, we implement the first method for incorporating grid rotations, wherein we use rotations which keep one axis fixed. Since the level sets of the solution are cubes, we will use rotations of $\beta = \pi/4$ in attempt to capture the edges. We will implement three rotations, alternately keeping the $x$-, $y$-, or $z$-axis fixed. The local upwind approximations are then \begin{equation} \label{eq:edgerotations} \begin{split}
\phi^{*,z}_{ijk} (a) &= \frac{1+ \frac{\abs{a_1 + a_2}}{\sqrt 2\Delta s_{xy}} \phi_{i+\text{sign}(a_1+a_2), j+\text{sign}(a_1+a_2),k} + \frac{\abs{a_2 - a_1}}{\sqrt 2\Delta s_{xy}} \phi_{i-\text{sign}(a_2-a_1), j+\text{sign}(a_2-a_1),k} + \frac{1}{\Delta z} \phi_{i,j,k+\text{sign}(a_3)} }{\frac{\abs{a_1 + a_2}}{\sqrt 2\Delta s_{xy}} + \frac{\abs{a_2 -a_1}}{\sqrt 2\Delta s_{xy}} + \frac 1 {\Delta z}}, \\
& \\
\phi^{*,y}_{ijk} (a) &= \frac{1+ \frac{\abs{a_1 + a_3}}{\sqrt 2\Delta s_{xz}} \phi_{i+\text{sign}(a_1+a_3), j,k +\text{sign}(a_1+a_3)} + \frac{1}{\Delta y} \phi_{i,j+\text{sign}(a_2),k} +  \frac{\abs{a_3 - a_1}}{\sqrt 2\Delta s_{xz}} \phi_{i-\text{sign}(a_3-a_1), j,k+\text{sign}(a_3-a_1)}  }{\frac{\abs{a_1 + a_3}}{\sqrt 2\Delta s_{xz}}  + \frac 1 {\Delta y} + \frac{\abs{a_3 -a_1}}{\sqrt 2\Delta s_{xz}}},\\
& \\ 
\phi^{*,x}_{ijk} (a) &= \frac{1+ \frac{1}{\Delta x} \phi_{i+\text{sign}(a_1),j,k} + \frac{\abs{a_2 + a_3}}{\sqrt 2\Delta s_{yz}} \phi_{i, j+\text{sign}(a_2+a_3),k +\text{sign}(a_2+a_3)}  +  \frac{\abs{a_3 - a_2}}{\sqrt 2\Delta s_{yz}} \phi_{i, j-\text{sign}(a_3-a_2),k+\text{sign}(a_3-a_2)}  }{ \frac{1}{\Delta x} +\frac{\abs{a_2 + a_3}}{\sqrt 2\Delta s_{yz}}  + \frac{\abs{a_3 -a_2}}{\sqrt 2\Delta s_{yz}}},
 \end{split}
\end{equation} where the superscript denotes the axis that is fixed, $\Delta s_{xy} = \sqrt{(\Delta x)^2 + (\Delta y)^2}$ and similarly for $\Delta s_{xz}$ and $\Delta s_{yz}$. The update rule for the iteration is then \begin{equation} \label{eq:updaterule3deik1edges} \phi^n_{ijk} = \min\bigg\{\phi^{n-1}_{ijk}, \,\,\,\min_{a}  \phi^*_{ijk}(a), \,\,\, \min_{a}  \phi^{*,x}_{ijk}(a), \,\,\, \min_{a}  \phi^{*,y}_{ijk}(a), \,\,\, \min_{a}  \phi^{*,z}_{ijk}(a)\bigg\}. \end{equation} The level set of the solution created using this update rule is seen in magenta in \cref{fig:lastsubfig2}. Note that while the edges are captured fairly sharply, the corners are still rounded off.   

In order to capture the corners sharply, we need to consider derivative approximations in the directions pointing toward the corners. It is a happy coincidence in 2D, that we can use a single rotation to capture all four corners of the square, since the vectors $(1,1)$ and $(1,-1)$ are orthogonal. In 3D, the vectors that point to alternate corners of the cube are no longer orthogonal; for example, $(1,1,1)$ is not orthogonal to $(1,-1,1)$.  Thus the rotation which captures the corners along the directions $(\pm1,\pm1,\pm1)$, will not capture any of the other corners.  Hence, if we want to capture all corners, we need to use four separate rotated approximations to the derivatives. 

We will describe the rotation that captures the corners in the directions of $(\hat \imath , \hat \jmath, \hat k) = (1, -1, 1)$, detailing every step along the way. To orient the $\overline x$-axis toward $(1,-1,1)$, we first rotate about the $z$-axis by an angle of $\beta = -\pi/4$, and then about the line $x=y$ by an angle of $\gamma = \arctan(1/\sqrt 2)$. The matrix that accomplishes this transformation is \begin{equation} \label{eq:1normeik3drotationmat}
U = \begin{pmatrix} \hphantom{-}1/\sqrt 3 & & \hphantom{-}1/\sqrt 2 & &  -1/\sqrt 6 \\ 
-1/\sqrt 3& & \hphantom{-}1/\sqrt 2& &\hphantom{-}1/\sqrt 6\\
 \hphantom{-}1/\sqrt 3 & & \hphantom{-}0 & & \hphantom{-}2/\sqrt 6
\end{pmatrix}
\end{equation} 
Following the computations above, the new grid directions are \begin{equation} \label{eq:3dDirectionsSpecificExample} \begin{split} &\text{positive } \overline x \text{ direction:} \,\ (1, -1, 1),\\ &\text{positive } \overline y \text{ direction:} \,\, (1, 1, 0), \\ &\text{positive } \overline z \text{ direction:} \,\, (-1, 1, 2), \end{split} \end{equation} and the rotated coefficient functions---which in this case depend only on $a$---are \begin{equation} 
\overline f_1(a) = \frac 1{\sqrt 3} (a_1 - a_2 + a_3), \,\,\,\,\,\,\,\,\,\, \overline f_2(a) = \frac{1}{\sqrt 2}(a_1 + a_2), \,\,\,\,\,\,\,\,\,\, \overline f_3(a) = \frac 1 {\sqrt 6} (-a_1 + a_2 + 2a_3), 
\end{equation} and thus the rotated equation is \begin{equation} \label{eq:rotated3deik1norm}-1 = \inf_{a_i \in \{\pm 1\} } \bigg\{ 
\left(\frac{a_1 - a_2 + a_3}{\sqrt 3}\right) \phi_{\overline x} +\left( \frac{a_1 + a_2}{\sqrt 2}\right) \phi_{\overline y} + \left(\frac{-a_1 + a_2 + 2a_3}{\sqrt 6}\right) \phi_{\overline z} \bigg\}.  \end{equation} At grid points $(x_i,y_j,z_k)$, the upwind derivative approximations are given by \begin{equation}\label{eq:3deikrotatedupwind}
\begin{split}
(\overline f_1(a)\phi_{\overline x})_{ijk} &= \frac{\abs{a_1 - a_2 + a_3}}{\sqrt 3} \frac{\phi_{i+\sign(a_1 - a_2 + a_3), j - \sign(a_1 - a_2 + a_3), k + \sign(a_1 - a_2 + a_3)} - \phi_{ijk}}{\Delta \overline x},\\
& \\
(\overline f_2(a)\phi_{\overline y})_{ijk} &= \frac{\abs{a_1 + a_2}}{\sqrt 2} \frac{\phi_{i+\sign(a_1 + a_2), j+ \sign(a_1 + a_2), k} - \phi_{ijk}}{\Delta \overline y}, \\
& \\
(\overline f_3(a)\phi_{\overline z})_{ijk} &= \frac{\abs{-a_1 + a_2 + 2a_3}}{\sqrt 6} \frac{\phi_{i-\sign(-a_1 + a_2 + 2a_3), j +\sign(-a_1 + a_2 + 2a_3), k +2\sign(-a_1 + a_2 + 2a_3)} - \phi_{ijk}}{\Delta \overline z},
\end{split}
\end{equation} where the new grid parameters are given by \begin{equation} \label{eq:3deik1normrotatedgridparameters}
\begin{split} 
\Delta \overline x &= \sqrt{(\Delta x)^2 + (\Delta y)^2 + (\Delta z)^2}, \\
\Delta \overline y &= \sqrt{(\Delta x)^2 + (\Delta y)^2}, \\
\Delta \overline z &= \sqrt{(\Delta x)^2 + (\Delta y)^2 + (2\Delta z)^2}.
\end{split}
\end{equation} Plugging these into the equation gives the upwind approximation to the equation at grid points: \begin{equation} \label{eq:upwindapproxtoequation3deik} \begin{split}
\phi^{*, +-+}_{ijk}(a) = &\bigg\{ 1 + \frac{\abs{a_1 - a_2 + a_3}}{\sqrt 3 \Delta x} \phi_{i+\sign(a_1 - a_2 + a_3), j - \sign(a_1 - a_2 + a_3), k + \sign(a_1 - a_2 + a_3)} \\ & \hspace{0.4cm}+ \frac{\abs{a_1 + a_2}}{\sqrt 2\Delta \overline y}\phi_{i+\sign(a_1 + a_2), j+ \sign(a_1 + a_2), k} \\ 
&\hspace{0.4cm} +\frac{\abs{-a_1 + a_2 + 2a_3}}{\sqrt 6\Delta \overline z} \phi_{i-\sign(-a_1 + a_2 + 2a_3), j +\sign(-a_1 + a_2 + 2a_3), k +2\sign(-a_1 + a_2 + 2a_3)}\bigg\}\\
&\hspace{0.2cm} / \bigg\{\frac{\abs{a_1 - a_2 + a_3}}{\sqrt 3 \Delta \overline x}+ \frac{\abs{a_1 + a_2}}{\sqrt 2\Delta \overline y} + \frac{\abs{-a_1 + a_2 + 2a_3}}{\sqrt 6\Delta \overline z}\bigg\},
\end{split}
\end{equation} where the superscript $+-+$ denotes the fact that the $\overline x$-axis points at the corners along the line parallel to $(1,-1,1)$. Finally, we can include this approximation, and iterate using the update rule \begin{equation}\label{eq:updateRule1Corner} \phi^n_{ijk} = \min\bigg\{\phi^{n-1}_{ijk}, \,\,\,\min_{a}  \phi^*_{ijk}(a), \,\,\, \min_{a}  \phi^{*,x}_{ijk}(a), \,\,\, \min_{a}  \phi^{*,y}_{ijk}(a), \,\,\, \min_{a}  \phi^{*,z}_{ijk}(a), \,\,\, \min_a \phi^{*,+-+}_{ijk}(a)\bigg\} \end{equation} In doing so, we will capture all of the edges of the level set fairly well, and perfectly capture the corners in the directions of $(1,-1,1)$. This is demonstrated by the cyan level set in \cref{fig:lastsubfig3}. Note that the remaining corners, along the directions $(1,1,1), (1,1, -1)$ and $(1,-1,-1)$ are still rounded, while the corners along $(1,-1,1)$ are sharp. 

Finally, if we want to perfectly capture all corners, we simply need to devise similar upwind approximations $\phi^{*,+++}_{ijk}(a),$ $\phi^{*,++-}_{ijk}(a),$  $\phi^{*,+--}_{ijk}(a)$, and include these using the update rule  \begin{equation}\label{eq:updateRuleAllCorners} \begin{split} \phi^n_{ijk} = \min\bigg\{\phi^{n-1}_{ijk}, \,\,\,&\min_{a}  \phi^*_{ijk}(a), \,\,\, \min_{a}  \phi^{*,x}_{ijk}(a), \,\,\, \min_{a}  \phi^{*,y}_{ijk}(a), \,\,\, \min_{a}  \phi^{*,z}_{ijk}(a), \,\,\,\\ &\min_a \phi^{*,+++}_{ijk}(a), \,\,\, \min_a \phi^{*,+-+}_{ijk}(a), \,\,\, \min_a \phi^{*,++-}_{ijk}(a) \min_a, \,\,\, \phi^{*,+--}_{ijk}(a)\bigg\}  \end{split}\end{equation} The level set of the solution resulting from this update rule is shown in yellow in \cref{fig:lastsubfig4}. In this case, the level set is (to machine precision) a perfect cube, with all corners and edges sharp. 

Lastly, \cref{tab:err3d} documents the maximal error in the numerical solution resolved using each update rule \eqref{eq:updaterule3deik1basic},\eqref{eq:updaterule3deik1edges},\eqref{eq:updateRule1Corner},\eqref{eq:updateRuleAllCorners}. As expected, including more approximations to the derivatives in additional directions only improves the accuracy. In the last trial, when we perfectly capture all edges and corners, the method is accurate to machine precision.  

\renewcommand{\arraystretch}{2}
\begin{table}[h!]
\centering
\begin{tabular}{l || c |  c |  c | c}
Update Rule &  \eqref{eq:updaterule3deik1basic} & \eqref{eq:updaterule3deik1edges} & \eqref{eq:updateRule1Corner} & \eqref{eq:updateRuleAllCorners} \\ [0.5ex]
\hline
Max Error & 1.0429e-01 &  4.2424e-02 & 3.9865e-02 &  1.4433e-14
\end{tabular}
\caption{Maximal error in our approximation to the solution of $\|\nabla \phi\|_1 = 1$ resulting from the update rules \eqref{eq:updaterule3deik1basic} (no rotated directions), \eqref{eq:updaterule3deik1edges} (rotations which keep one axis fixed), \eqref{eq:updateRule1Corner} (an additional rotation to resolve the corners along $(1,-1,1)$), and \eqref{eq:updateRuleAllCorners} (rotations to capture all the corners). Each approximation was computed on a $201 \times 201 \times 201$ grid. Notice that when all corners and edges are accounted for, the solution is accurate to machine precision.}
\label{tab:err3d}
\end{table}

\end{document}